                      \def\version{19 June, 2009}                           %

\documentclass[reqno,12pt]{amsart} 

\usepackage{amsmath} 
\usepackage{amssymb} 
 
\newenvironment{proofsect}[1] 
{\vskip0.1cm\noindent{\bf #1.}\hskip0.5cm}

\def\ti{\to\infty}

\def\d{\delta} 
\def\e{\varepsilon} 

\def\n{\nu} 
 
\newfam\Bbbfam 
\font\tenBbb=msbm10 
\font\sevenBbb=msbm7 
\font\fiveBbb=msbm5 
\textfont\Bbbfam=\tenBbb 
\scriptfont\Bbbfam=\sevenBbb 
\scriptscriptfont\Bbbfam=\fiveBbb

\newcommand{\R}     {\mathbb{R}} 
 
\newcommand{\N}     {\mathbb{N}}

\def\1{{\mathchoice {1\mskip-4mu\mathrm l}      
{1\mskip-4mu\mathrm l} 
{1\mskip-4.5mu\mathrm l} {1\mskip-5mu\mathrm l}}} 
\newcommand{\ssup}[1] {{{\scriptscriptstyle{({#1}})}}} 
\def\comment#1{} 
\newtheoremstyle{thm}{2ex}{2ex}{\itshape\rmfamily}{} 
{\bfseries\rmfamily}{}{1.7ex}{} 
 
\newtheoremstyle{rem}{1.3ex}{1.3ex}{\rmfamily}{} 
{\itshape\rmfamily}{}{1.5ex}{}

\newtheorem{theorem}{Theorem}[section] 
\newtheorem{Lemma}[theorem]{Lemma} 
\newtheorem{proposition}[theorem] {Proposition} 
\newtheorem{cor}[theorem]  {Corollary} 
 
\newtheorem{definition}[theorem] {Definition}

\theoremstyle{definition}

\renewcommand{\section}{\secdef\sct\sect} 
\newcommand{\sct}[2][default]{\refstepcounter{section} 
\vspace{0.8cm} 
\setcounter{equation}{0} 
\centerline{ 
\large\scshape \arabic{section}.\ #1} 
\vspace{0.2cm}} 
\newcommand{\sect}[1]{ 
\vspace{0.8cm} 
\centerline{\large\scshape #1} 
\vspace{0.2cm}} 
 
\renewcommand{\subsection}{\secdef \subsct\sbsect} 
\newcommand{\subsct}[2][default]{\refstepcounter{subsection} 
\nopagebreak 
\vspace{0.5\baselineskip} 
{\flushleft\bf \arabic{section}.\arabic{subsection}~\bf #1  } 
\nopagebreak} 
\newcommand{\sbsect}[1]{\vspace{0.1cm}\noindent 
{\bf #1}\vspace{0.1cm}}

\renewcommand{\subsubsection}{%
\secdef \subsubsect\sbsbsect} 
\newcommand{\subsubsect}[2][default]{%
\refstepcounter{subsubsection} 
\nopagebreak 
\vspace{0.1\baselineskip} 
\nopagebreak 
{\flushleft 
\sffamily\slshape 
\arabic{section}.\arabic{subsection}.\arabic{subsubsection} 
\ %
\sffamily #1\/.}\ } 
\newcommand{\sbsbsect}[1]{\vspace{0.1cm}\noindent 
{\bf #1}\ } 
 
\renewcommand{\d}{{\rm d}} 
 
\newcommand{\eps}{\varepsilon} 
 
\newcommand{\Sym}{\mathfrak{S}}

\newcommand{\Dcal}   {{\mathcal D }} 
 
\newcommand{\Fcal}   {{\mathcal F }} 
\newcommand{\Gcal}   {{\mathcal G }} 
\newcommand{\Hcal}   {{\mathcal H }} 
 
\newcommand{\Jcal}   {{\mathcal J }} 
 
\newcommand{\Lcal}   {{\mathcal L }}

\newcommand{\Tcal}   {{\mathcal T }}

\renewcommand{\e}   {{\operatorname e }}
 
\begin{document}
\title[Limit theorems for VRJP on trees] {Limit  theorems  for vertex-reinforced jump processes on regular trees}
\author[Andrea Collevecchio]{}

\maketitle
\begin{center}
\thispagestyle{empty} 
\vspace{0.2cm} 
 
\centerline {\sc By  Andrea Collevecchio\footnote{Dipartimento di Matematica applicata, Universit\`a Ca' Foscari -- Venice, Italy.  {\tt collevec@unive.it}}}

\vspace{0.4cm}

\centerline{\small(\version)}
\vspace{.5cm}
\end{center}
\[\]
\begin{abstract}
Consider a vertex-reinforced jump process  defined on a regular tree, where each vertex has exactly $b$ children, with $ b \ge 3$. We prove the strong law of large numbers  and the central limit theorem for the distance of the process from the root.  Notice that it is still unknown if vertex-reinforced jump process  is transient on the binary tree.  
\end{abstract}
\newpage
\section{Introduction}
\label{intro}

\noindent  Let $\Dcal$ be any graph with the property that each vertex is the end point of only a finite number of edges. Denote by  Vert$(\Dcal)$ the set of vertices of $\Dcal$. The following, together with the vertex occupied at time $0$ and the set of positive numbers $\{a_{\nu} \colon \nu \in \mbox{Vert}(\Dcal) \}$,  defines a right-continuous process  $\mathbf{X} = \{ X_s,\; s \ge 0\} $.  This process  takes as values the vertices of $\Dcal$  and jumps only to nearest neighbors, i.e. vertices one edge away from the occupied one. Given $X_s, \; 0 \le s \le t $, and  $ \{ X_t = x \}$, the conditional probability that,  in the interval $(t,t+\d t)$, the process jumps to the nearest neighbor $y$ of $x$  is $L(y,t) \d t$, with
\[ L(y, t) := a_{y} + \int_0^t \1_{\{X_s = y \} } \d s, \qquad a_{y} >0,\]
where $\1_{A}$ stands for the indicator function of the set $A$.
 The positive numbers $\{a_{\nu} \colon \nu \in \mbox{Vert}(\Dcal)\}$ are called  initial weights, and we suppose $ a_{\nu} \equiv 1$, unless specified otherwise. Such a process is said to be a  Vertex Reinforced Jump Process (VRJP) on $\Dcal$. 

 Consider VRJP defined on the integers, which starts from 0. With probability $1/2$ it will jump   either to 1 or $-1$. The time of the first jump is  an exponential random variable with mean $1/2$, and is independent on the direction of the jump. Suppose the walk jumps towards 1 at time $z$. Given this, it will wait at 1 an exponential amount of time with mean $1/(2+z)$. Independently of this time, the jump will be towards $0$ with probability $(1+z)/(2+z)$.
 
 In this paper we define a process to be  recurrent  if it visits each vertex infinitely many times  a.s., and to be  transient otherwise. VRJP was introduced by Wendelin Werner, and its properties were first studied by Davis and Volkov (see   \cite{DV2002} and \cite{DV2004}).   This reinforced walk defined on the integer lattice is studied in \cite{DV2002} where recurrence is proved. For fixed $b \in \N :=\{1,2, \ldots\}$, the   $b$-ary tree,  which we denote by $\Gcal_{b}$, is the infinite tree where each vertex has $b+1$ neighbors with the exception of a single vertex, called the root and designated by $\rho$, that
is connected to $b$ vertices.  In \cite{DV2004} is shown that VRJP on the $b$-ary tree is transient if  $ b \ge 4$. The case $b=3$ was dealt in \cite{C2006b}, where it was proved that the process is still transient. The case $b=2$ is still open. 

Another  process which reinforces the vertices, the so called  Vertex-Reinforced Random Walk (VRRW), shows a completely different behaviour. VRRW was introduced by Pemantle (see \cite{P1992}). Pemantle and Volkov (see \cite{PV99}) proved  that this process, defined on the integers, gets stuck in at most five points. Tarr\`es (see \cite{T2004}) proved that it gets stuck in exactly 5 points. Volkov (in \cite{V2001}) studied this process on arbitrary trees.  

The reader can find in \cite{P2007}  a survey on reinforced processes. 
In particular, we would like to mention that little is known regarding the behaviour of these processes on infinite graphs with loops.   Merkl and Rolles (see \cite{MR2008}) studied the recurrence of  the original reinforced random walk, the so-called linearly bond-reinforced random walk,  on two-dimensional graphs. Sellke (see \cite{S2006}) proved than once-reinforced random walk is recurrent on the ladder.
  
We define the distance between two vertices as the number of edges in the unique self-avoiding path connecting them. For any vertex $\nu$, denote by $|\nu|$  its distance from the root. Level $i$ is  the set of vertices $\n$ such that $|\nu|=i$. The main result of this paper is the following.


\begin{theorem}\label{stronglaw} 
Let $\mathbf{X}$ be VRJP on $\Gcal_{b}$, with $ b \ge 3$. There exist constants $  K^{\ssup 1}_b  \in (0, \infty)$ and $ K^{\ssup 2}_b \in [0, \infty)$ such that   
\begin{eqnarray}
&&\label{slln} \lim_{t \ti }\frac{|X_t|}{t} = K^{\ssup 1}_b \qquad \mbox{a.s.}, \\
&&\frac{ |X_t| - K^{\ssup 1}_b t}{ \sqrt{t}} \Longrightarrow  \mbox{Normal$(0, K^{\ssup 2}_b)$},
\end{eqnarray}
where we took the limit as $t \ti$, $ \Rightarrow$ stands for weak convergence and Normal$(0,0)$ stands for the Dirac mass at $0$.
\end{theorem}

\noindent  Durrett, Kesten and Limic have proved in \cite{DKL2002} an analogous result  for a bond-reinforced random walk, called one-time bond-reinforced random walk, on $\Gcal_{b}$, $ b \ge 2$. To prove this, they break the path into independent identically distributed  blocks, using the classical method of cut points. We also use this approach.  Our implementation of the cut point method is a strong  improvement of the one used in \cite{C2006} to prove the strong law of large numbers for the original reinforced random walk, the so-called linearly bond-reinforced random walk, on $\Gcal_{b}$, with $b \ge 70$.   Aid\'ekon, in \cite{A2008} gives a  sharp criteria for random walk in a random environment, defined on Galton-Watson tree, to have positive speed. He proves the strong law of large numbers for linearly bond-reinforced random walk on $\Gcal_{b}$, with $b \ge 2$.


\section{Preliminary definitions and properties}
From now on, we consider VRJP $\mathbf{X}$ defined on the regular tree $\Gcal_{b}$, with $b \ge 3$. For $\nu \neq \rho$, define $\mbox{par}(\nu)$, called the parent of $\nu$,  to be the unique vertex  at level  $|\nu|-1$ connected to $\nu$.  A vertex $\nu_0$ is a child of $\nu$ if $\nu = \mbox{par}(\nu_0)$. We say that a vertex $ \nu_{0}$ is a  descendant of the vertex $\nu$ if the latter lies on the unique self-avoiding path connecting $\nu_{0}$ to $\rho$, and $ \nu_{0} \neq \nu$. In this case, $\nu$ is said to be an ancestor of $\nu_{0}$.
  For any vertex $\mu$, let $\Lambda_{\mu}$  be the subtree   consisting of $\mu$, its descendants and the edges connecting them, i.e. the subtree rooted at $\mu$.  Define
\begin{eqnarray*}
 T_i &:=& \inf \{ t \ge 0 \colon |X_t| = i \}.
\end{eqnarray*} 

We give the so-called Poisson construction of VRJP on a graph $\Dcal$ (see \cite{S1994}). For each ordered pair of neighbors $ (u,v)$ assign a Poisson process $P(u,v)$ of rate 1, the processes being independent. Call $h_{i}(u,v)$, with $ i \ge 1$, the inter-arrival times of $ P(u,v)$ and let $ \xi_{1} := \inf \{ t \ge 0 \colon X_{t} = u\}$. The first jump after $\xi_{1}$ is at time $ c_{1} := \xi_{1} + \min_{v} h_{1}(u,v) \big(L(v, \xi_{1})\big)^{-1}$, where the minimum is taken over the set of neighbors of $u$. The jump is towards the neighbor $v$ for which that minimum is attained. Suppose we defined $\{ (\xi_{j}, c_{j}), 1 \le j \le i-1 \}$, and let 
\begin{eqnarray*}
 &&\xi_{i} := \inf \big\{ t > c_{i-1} \colon X_{t} = u \big\}, \mbox{ and} \\
 && j_{v} -1 = j_{u,v} -1 := \mbox{ number of times $\mathbf{X}$ jumped from $u$ to $v$ by time $\xi_{i}$}.
 \end{eqnarray*}
The first jump after $\xi_{i}$ happens at time $ c_{i}:= \xi_{i} + \min_{v} h_{j_{v}}(u,v) \big(L(v, \xi_{i})\big)^{-1}$, and the jump is towards the neighbor $v$ which attains that minimum.
\begin{definition}
A vertex $\mu$, with $|\mu| \ge 2$, is {\bf good} if it satisfies the following
\begin{equation}\label{good}
h_{1}(\mu_{0}, \mu) < \frac{h_{1}\big(\mu_{0}, \mbox{\rm par}(\mu_{0})\big)}{1 + h_{1}\big( \mbox{\rm par}(\mu_{0}), \mu_{0} \big)}\qquad \mbox{ where } \mu_{0}= \mbox{\rm par}(\mu).
\end{equation}
\end{definition}
By virtue of  our construction of VRJP, \eqref{good} can be interpreted as follows. When the process $\mathbf{X}$ visits the vertex $\mu_{0}$ for the first time, if this ever happens, the weight at its parent is exactly $1 + h_{1}\big( \mbox{\rm par}(\mu_{0}), \mu_{0}\big)$ while the weight at $\mu$ is $1$. Hence condition   \eqref{good} implies that  when the process visits  $\mu_{0}$ (if this ever happens)   then it will visit $\mu$ before it returns to par($\mu_{0}$), if this ever happens.

The next Lemma gives bounds for the probability that VRJP returns to the root after the first jump.
\begin{Lemma}\label{return}
 Let
\[\alpha_b := \mathbb{P} \bigl( X_t = \rho \mbox{ for some $ t\ge T_1$}\bigr), \] 
and let $\beta_b$ be the smallest among the positive  solutions of the equation
\begin{equation}
x= \sum_{k=0}^{b} x^{k} p_{k},
\end{equation}
where, for $k \in \{0,1, \ldots, b\}$, 
\begin{equation}\label{pikappa}
p_{k}:= \sum_{j=0}^{k} { b \choose k}  {k \choose j} (-1)^{j}  \int_0^{\infty}  \frac{1 +z}{j+b-k+1+z} \e^{-z}  \d z.\qquad
\end{equation}
 We have 
\begin{equation} \label{goodret}
\int_{0}^{\infty} \frac{1+z}{b+1+z}b\e^{-bz}\d z \le \alpha_{b} \le \beta_{b}. 
\end{equation} 
\end{Lemma}
\begin{proofsect}{Proof} First we prove the lower bound in \eqref{goodret}.  The left-hand side of  this inequality is the probability that the process returns to the root with exactly two jumps. To see this, notice that $L(\rho, T_{1})$ is equal $ 1 +  \min_{\nu \colon |\nu|=1} h_{1}(\rho, \nu)$. Hence $T_{1}=L(\rho, T_{1})-1$ 
 is distributed like an exponential with mean $1/b$.  Given that $T_{1}=z$, the probability that the second  jump is from $X_{T_{1}}$ to $\rho$ is equal to $(1+z)/(b+1+z)$.  Hence the probability that the process returns to the root with exactly two jumps is
\begin{equation*}
 \int_{0}^{\infty} \frac{1+z}{b+1+z}b\e^{-bz}  \d z.
\end{equation*} 
As for the upper bound in \eqref{goodret} we reason as follows. We give an upper bound for the probability that  there exists  an infinite random tree which is composed only of good vertices and which has root at one of the children of $X_{T_{1}}$. If this event holds, then the process does not return to the root after time $T_{1}$ (see the proof of Theorem 3 in \cite{C2006b}). We prove that a  particular cluster of good vertices is  stochastically larger than  a  branching process which is supercritical. We introduce the following color scheme. The only vertex at level 1 to be {\it green} is $X_{T_{1}}$. A vertex $\nu$, with $|\nu|\ge 2$,  is {\it green} if and only if it is good and its parent is green. All the other vertices are uncolored.
 Fix  a vertex $\mu$. Let $ C$ be any event in 
\begin{equation}\label{outsigma}
 \Hcal_{\mu}:= \sigma( h_i(\eta_0, \eta_1): i \ge 1,   \mbox{ with $ \eta_0 \sim \eta_1 $ and  both } \eta_0 \mbox{ and } \eta_1 \notin \Lambda_{\mu} ),
 \end{equation}
  that is  the $ \sigma$-algebra that contains the information about  $X_t$   observed  outside $ \Lambda_{\mu}$. Next we show that given $C \cap \{ \mu $ is green$\}$, the distribution of $ h_1(\mbox{par}(\mu), \mu)$ is stochastically dominated by an exponential(1). To see this, first notice that $ h_1(\mbox{par}(\mu), \mu)$ is independent of $ C$. Let $ D:= \{ \mbox{par}(\mu) \mbox{ is green}\} \in \Hcal_{\mu}$ and set
 \begin{equation}\label{stucapdiw}
 W := \frac{h_{1}\big(\mu_{0}, \mbox{\rm par}(\mu_{0})\big)}{1 + h_{1}\big( \mbox{\rm par}(\mu_{0}), \mu_{0} \big)}\qquad \mbox{ where } \mu_{0}= \mbox{\rm par}(\mu).
\end{equation}
The random variable $W$ is independent of 
  $ h_1(\mbox{par}(\mu), \mu) $ and is absolutely continuous with respect the Lebesgue measure.  By the definition   of good vertices we have   
\[  \{ \mu \mbox{ is green}\} = \{ h_1(\mbox{par}(\mu), \mu) < W \} \cap D.\]
\noindent Denote by $f_{W}$ the  conditional density of $W$ given $D \cap C\cap \{ h_1(\mbox{par}(\mu), \mu) < W \}$. We have
\begin{equation}\label{cluster2}
\begin{aligned}
&\mathbb{P} \Bigl(  h_1(\mbox{par}(\mu), \mu) \ge x    \;  \big{|} \;  \{\mbox{$\mu$ is green}\}\cap C \Bigr) \\
=\quad  &\mathbb{P} \Bigl(  h_1(\mbox{par}(\mu), \mu) \ge x    \; \big{|} \; \{ h_1(\mbox{par}(\mu), \mu) < W \}   \cap C \cap D \Bigr) \\ 
=\quad &\int_{0}^{\infty} \mathbb{P} \Bigl(  h_1(\mbox{par}(\mu), \mu) \ge x    \;  \big{|} \; \{ h_1(\mbox{par}(\mu), \mu) < w \}   \cap C \cap D  \cap \{ W = w\} \Bigr) f_{W}(w) \d w 
\end{aligned}
\end{equation}
Using the facts that $h_{1}($par$(\mu),\mu)$ is independent of $W, C$ and $D$ and $$\mathbb{P}( h_1(\mbox{par}(\mu), \mu) \ge x \;|\;  h_1(\mbox{par}(\mu), \mu) < w) \le \mathbb{P}( h_1(\mbox{par}(\mu), \mu) \ge x ),$$
 we get that the  expression in \eqref{cluster2} is 
 less or equal to $\mathbb{P} \Bigl(  h_1(\mbox{par}(\mu), \mu) \ge x   \Bigr)$.
Summarising 
\begin{equation}\label{sumclust}
 \mathbb{P}\Big( h_1(\mbox{par}(\mu), \mu)\ge x \;|\; \{\mu \mbox{ is green}\}\cap C\Big) \ge  \mathbb{P} \Bigl(  h_1(\mbox{par}(\mu), \mu) \ge x   \Bigr).
\end{equation}
The inequality \eqref{cluster2} implies that  if $\mu_{1}$ is a child of $\mu$ and $C \in \Hcal_{\mu}$ we have
\begin{equation}
\mathbb{P} \Bigl( \mbox{$\mu_{1}$ is green}    \;  \mid \;  \{\mbox{$\mu$ is green}\}\cap C \Bigr) \ge 
\mathbb{P} \Bigl( \mbox{$\mu_{1}$ is green} \Bigr).
\end{equation}
To see this, it is enough to integrate over the value of $h_1(\mbox{par}(\mu), \mu)$ and use the fact that, conditionally on $h_1(\mbox{par}(\mu), \mu)$, the events $\{ \mbox{$\mu_{1}$ is green} \}$ and  $\{\mbox{$\mu$ is green}\}\cap C$ are independent. The probability that $\mu_{1}$ is good conditionally on   $\{h_1(\mbox{par}(\mu), \mu)=x\}$ is a non-increasing function of $x$, while the  distribution of $h_1(\mbox{par}(\mu), \mu)$ is stochastically smaller than the   conditional distribution of $h_1(\mbox{par}(\mu), \mu)$ given $\{\mbox{$\mu$ is green}\}\cap C$, as shown in \eqref{sumclust}.

Hence the cluster of green vertices is stochastically larger than a Galton--Watson tree where each vertex has $k$ offspring,  $k \in \{0,1, \ldots, b\}$, with probability $p_k$ defined in \eqref{pikappa}. To see this, fix a vertex $\mu$ and let $\mu_{i}$, with $i \in \{0,1, \ldots, b\}$ be its children.  It is enough to realize that $p_{k}$ is  the probability that exactly $k$ of the $h_{1}(\mu, \mu_{i})$, with $i \in \{0,1, \ldots, b\}$, are smaller than $\big(1 + h_{1}(\mbox{par}(\mu),\mu) \big)^{-1}h_{1}\big(\mu, \mbox{par}(\mu)\big)$. As the random variables $h_{1}(\mu, \mu_{i}), h_{1}\big(\mu, \mbox{par}(\mu)\big)$ and $h_{1}(\mbox{par}(\mu),\mu)$ are independent exponentials with parameter one, we have 
\begin{equation}
\begin{aligned}
p_{k} &= {b \choose k} \int_{0}^{\infty} \int_{0}^{\infty} \mathbb{P}\big(h_{1}(\mu_{0}, \mu)< \frac{y}{1 +z} \big)^{k} \mathbb{P}\big(h_{1}(\mu_{0}, \mu)\ge \frac{y}{1 +z} \big)^{b-k} \e^{-y} \e^{-z} \d y \, \d z\\
&= {b \choose k} \int_{0}^{\infty} \int_{0}^{\infty} \big(1- \e^{- \frac{y}{1+z}}\big)^{k} \e^{- \frac{y}{1+z} (b-k)} \e^{-y} \e^{-z} \d y \, \d z \\
&= \sum_{j=0}^{k} \int_{0}^{\infty} \int_{0}^{\infty} {b \choose k} {k \choose j} (-1)^{j} \e^{-y(j+b -k +1+z)/(1+z)} \e^{-z} \d y \, \d z\\
&= \sum_{j=0}^{k} {b \choose k} {k \choose j} (-1)^{j} \int_{0}^{\infty} \frac{1+z}{j+b -k +1 +z} \e^{-z} \d z.
\end{aligned}
\end{equation}

From the basic theory of branching processes we know that the probability that this Galton--Watson tree is finite (i.e. extinction) equals the smallest positive solution of the equation 
\begin{equation}\label{cluster4}
x -\sum_{k=0}^b x^k p_k =0.
\end{equation}
The proof of \eqref{goodret} follows from the fact that  $1-\beta_{b } \le 1-\alpha_{b}$.
This latter inequality is a consequence of the fact that the cluster  of green vertices  is stochastically larger than the Galton-Watson tree, hence its probability of non-extinction is not smaller. As $b \ge 3$, the Galton-Watson tree is supercritical (see \cite{C2006b}),hence $\beta_{b}<1$.
\qed           \end{proofsect}

For example,  if we consider VRJP on $\Gcal_{3}$, Lemma~\ref{return} yields 
  $$
  0.3809 \le \alpha_{3} \le 0.8545.
  $$

\begin{definition}    
 Level $j \ge 1 $ is a  {\bf cut level}     if  the first jump after $T_{j}$ is towards level $j+1$, and after time $T_{j+1}$   the process never goes back to $X_{T_j}$, and 
$$  L( X_{T_j}, \infty ) < 2 \quad \mbox{and} \quad L(\mbox{\rm par}( X_{T_j}), \infty ) < 2. $$ 
 Define $l_1$ to be the  cut level with minimum distance from the root, and for $i >1$,
\[ l_i := \min \{ j > l_{i-1}\colon j \mbox{ is a cut level}\}. \]
\item Define the $i$-th {\bf cut time} to be $\tau_i  := T_{l_i} $. Notice that $ l_i = | X_{\tau_i}|$.
\end{definition}



\section{ $l_1$  has an exponential tail} 
For any vertex $\nu \in \mbox{Vert}(\Gcal_{b})$,  we define fc($\nu$), which stands for {\bf first child} of $\nu$, to be the (a.s.) unique vertex connected to $\nu$ satisfying
 \begin{equation}\label{firstc}
 h_{1}(\nu, \mbox{fc}(\nu)) = \min \big\{ h_{1}(\nu, \mu) \colon \mbox{par}(\mu) = \nu \big\}.
 \end{equation}
For definiteness, the root $\rho$ is not a first child. Notice that condition \eqref{firstc} does not imply that the vertex fc($\nu$) is visited by the process. If $\mathbf{X}$ visits it, then it is the first among the children of $\nu$ to be visited.

For any pair of distributions $f$ and $g$, denote by $f \,\overline{*}\,g$ the distribution of $\sum^{V}_{ k = 1} M_{k}$, where 
\begin{itemize}
\item $V$ has distribution  $f$, and 
\item $\{ M_{k},\, k \in  \N \}$ is a sequence of i.i.d random variables, independent of $V$, each with distribution $g$.
\end{itemize}
 Recall  the definition of $p_{i}$, $i \in \{0, \ldots, b\}$, given in \eqref{pikappa}.  Denote by $\mathbf{p}^{\ssup 1}$ the distribution which assigns to $i \in \{0, \ldots, b\}$ probability $p_{i}$. Define, by recursion, $\mathbf{p}^{\ssup j} :=  \mathbf{p}^{\ssup {j-1}} \,\overline{*}\, \mathbf{p}^{\ssup 1}$, with $j \ge 2$. The distribution $\mathbf{p}^{\ssup j}$ describes the number of elements, at time $j$, in a population which evolves like a branching process generated by one ancestor and  with offspring distribution $\mathbf{p}^{\ssup 1}$. If we let
\[ 
m := \sum_{j=1}^{b} j p_{j},
\]
then the mean of $\mathbf{p}^{\ssup j}$ is $m^{j}$. The probability that a given vertex $\mu$ is good is, by definition, 
$$
\mathbb{P} \Big(h_{1}(\mu_{0}, \mu) < \frac{h_{1}\big(\mu_{0}, \mbox{\rm par}(\mu_{0})\big)}{1 + h_{1}\big( \mbox{\rm par}(\mu_{0}), \mu_{0} \big)}\Big) \qquad \mbox{ where } \mu_{0}= \mbox{\rm par}(\mu).
 $$
As the $ h_{1}\big({\rm par}(\mu_{0}), \mu_{0} \big)$ is exponential with parameter 1, conditioning on its value and using independence between different Poisson processes, we have that the probability above equals
\begin{equation}
\begin{aligned}
\mathbb{P}\Big(h_{1}(\mu_{0}, \mu) < \frac1{1 + z}h_{1}\big(\mu_{0}, \mbox{\rm par}(\mu_{0})\big)\Big) \e^{-z} \d z=
 \int_{0}^{\infty} \frac 1{2+z} \e^{-z} \d z= 0.36133\ldots.
\end{aligned}
\end{equation}
Hence
$$ m = b \cdot 0.36133 > 1,$$
because we assumed $b \ge 3$.

Let $q_{0} = p_{0} + p_{1},\, $ and for $k \in \{1,2, \ldots,b-1\}$ set    $ q_{k} = p_{k+1}$. Set $\mathbf{q}$ to be the distribution which assigns to $i \in \{0, \ldots, b-1\}$ probability $q_{i}$.  For $j \ge 2$, let
$\mathbf{q}^{\ssup j} := \mathbf{p}^{\ssup {j-1}} \,\overline{*}\, \mathbf{q}$. Denote by $q^{\ssup j}_{i}$ the probability that the distribution $ \mathbf{q}^{\ssup j}$ assigns to $i \in \{ 0, \ldots, (b-1)b^{j-1}\}$. The mean of $\mathbf{q}^{\ssup j}$ is $m^{j-1}(m-1)$. From now on, $\zeta$ denotes the smallest positive integer  in $\{2,3, \ldots,\}$   such that 
\begin{equation}\label{offwthfc}
m^{\zeta-1}(m-1) > 1.
\end{equation}

Next we want to define a sequence of events which are independent and which are closely related to the event that a given level is a cut level. 
For any vertex $\nu $ of $\Gcal_{b}$  let $\Theta_{\nu}$ be the set of vertices  $\mu$ such that 
\begin{itemize}
\item $\mu$ is a descendant of $\nu$, 
\item the difference $|\mu|$ - $|\nu|$ is a multiple of $\zeta$,
\item $\mu$ is a first child. 
\end{itemize}
 By subtree rooted at $\nu$ we mean   a subtree of $\Lambda_{\nu}$ that contains $\nu$. Set $\widetilde{\nu} =$ fc$(\nu)$ and   let
 \begin{equation}\label{sub1}
 \begin{aligned}
 A(\nu)&:= \big\{ \exists \mbox{ an infinite subtree of $\Gcal_{b}$ root at a child of $\widetilde{\nu}$, which is composed only by}\\
 &\qquad \mbox{good vertices and which contains none of the  vertices in $\Theta_{\nu}$}\}
 \end{aligned}
 \end{equation}
 For $ i \in \N$, let $A_{i} := A\big(X_{T_{i}}\big)$. Notice that if the process reaches the first child of $\nu$ and if $A(\nu)$ holds, then the process will never return to $\nu$. Hence if $A_{i}$ holds, and if $X_{T_{i+1}} = X_{T_{i}}+1$, then $i$ is a cut level, provided that the total weights at $X_{T_{i}}$ and its parent are less than 2.
 \begin{proposition}\label{aind}
 The events $A_{i \zeta}$, with $i \in \N$,  are independent.
 \end{proposition}
 \begin{proofsect}{Proof}
 We recall  that $\zeta \ge 2$. We proceed by backward recursion and show that the events $A_{i \zeta}$ depend on disjoint Poisson processes collections. Choose  integers $ 0 < i_{1}<i_{2}< \ldots<i_{k}$, with $i_{j} \in \zeta \N := \{ \zeta, 2 \zeta, 3 \zeta, \ldots\}$ for all $j \in \{1, 2, \ldots, k\}$.  It is enough to prove that
\begin{equation}\label{indep1}
\mathbb{P} \Big( \bigcap_{j=1}^{k} A_{i_{j}} \Big) = \prod_{j=1}^{k} \mathbb{P} \big( A_{i_{j}}\big).
\end{equation} 
Fix a vertex $ \nu$ at level $i_{k}$. 
The set $ A(\nu)$ belongs to the sigma-algebra generated by $\big\{  P(u,w) \colon \; u,w \in \mbox{Vert}(\Lambda_{\nu}) \big\}$. On the other hand, the set $\bigcap_{j=1}^{k-1} A_{i_{j}} \cap \{X_{T_{i_{k}}} = \nu\}$ belongs to $\big\{  P(u,w) \colon \; u \notin \mbox{Vert}(\Lambda_{\nu}) \big\}$. As the two events belong to disjoint collections of independent Poisson processes, they are independent. As $\mathbb{P}(A(\nu)) =  \mathbb{P}(A(\rho))$, we have
\begin{equation}
\begin{aligned}
& \mathbb{P} \Big( A_{i_{k}} \cap \bigcap_{j=1}^{k-1} A_{i_{j}} \Big) = \sum_{\nu \colon |\nu|=i_{k}} \mathbb{P}\Big(A_{i_{k}}  \cap  \bigcap_{j=1}^{k-1} A_{i_{j}} \cap \{X_{T_{i_{k}}} = \nu\}\Big)\\
&= \sum_{\nu \colon |\nu|=i_{k}} \mathbb{P}\Big(A(\nu)  \cap  \bigcap_{j=1}^{k-1} A_{i_{j}} \cap \{X_{T_{i_{k}}} = \nu\}\Big)= \sum_{\nu \colon |\nu|=i_{k}} \mathbb{P}\big(A(\nu) \big)\mathbb{P} \Big(  \bigcap_{j=1}^{k-1} A_{i_{j}} \cap \{X_{T_{i_{k}}} = \nu\} \Big) \\ 
&= \mathbb{P}\big(A(\rho) \big) \sum_{\nu \colon |\nu|=i_{k}} \mathbb{P} \Big(  \bigcap_{j=1}^{k-1} A_{i_{j}} \cap \{X_{T_{i_{k}}} = \nu\} \Big)=\mathbb{P}\big(A(\rho) \big) \mathbb{P} \Big(  \bigcap_{j=1}^{k-1} A_{i_{j}} \Big).
\end{aligned}
\end{equation}
The events $A(\nu)$ and $\{X_{T_{i_{k}}} = \nu\}$ are independent, and by virtue of the self-similarity property of the regular tree we get $\mathbb{P}\big(A(\rho) \big) = \mathbb{P}\big(A_{i_{k}} \big)$. Hence
\begin{equation}\label{recstep}
\mathbb{P} \Big( A_{i_{k}} \cap \bigcap_{j=1}^{k-1} A_{i_{j}} \Big)= \mathbb{P}\big(A_{i_{k}} \big) \mathbb{P} \Big(  \bigcap_{j=1}^{k-1} A_{i_{j}} \Big).
\end{equation}
Reiterating \eqref{recstep} we get \eqref{indep1}.
\qed
\end{proofsect}


 \begin{Lemma}\label{return1}
 Define $\gamma_{b}$ to be the smallest positive solution of the equation
 \begin{equation}\label{root}
 x = \sum_{k=0}^{b-1} x^{k}q^{\ssup \zeta}_{k},
 \end{equation}
 where $\zeta$ and $(q_{k}^{\ssup n})$ have been defined at the beginning of this section.
 We have 
 \begin{equation}\label{ret1}
 \mathbb{P}(A_{i}) \ge 1-\gamma_{b}>0, \qquad \forall i \in \N.
 \end{equation}
 \end{Lemma}

 \begin{proofsect}{Proof}
 Fix $i \in \N$ and let $\nu^{*}= X_{T_{i}}$.
   We adopt the following color scheme. The vertex fc$\big(X_{T_{i}}\big)$ is colored {\it blue}.  A descendant  $\mu$ of $\nu^{*}$ is colored blue if it is good,  its parent is {\it blue},  and either
\begin{itemize}
\item $|\mu|- |\nu^{*}|$ is not a multiple of $\zeta$,  or
\item $\frac 1\zeta \big( |\mu|- |\nu^{*}|\big) \in \N$ and $\mu$ is not a first child.
\end{itemize}
Vertices which are not descendants of $\nu^{*}$ are not colored.
 Following the reasoning given in the proof of Lemma \ref{return}, we can conclude that the number of blue vertices at levels $ |\nu^{*}|~+~ j\zeta$, with $j \ge 1$, is  stochastically larger than the number of individuals in a population which evolves like a branching process with offspring distribution  $\mathbf{q}^{\ssup \zeta}$, introduced at the beginning of this section.  Again, from the basic theory of branching processes we know that the probability that this tree is finite equals the smallest positive solution of the equation \eqref{root}. By virtue of \eqref{offwthfc} we have that $\gamma_{b}<1$.
 \qed           \end{proofsect}
 The proof of the following Lemma can be found in \cite{DZ}  pages 26-27 and 35. 
\begin{Lemma}\label{bind}
Suppose $U_n $ is  Bin$(n,p)$. For $  x \in (0,1)$ consider the entropy
$$H(x\,|\, p) := x \ln \frac{x}{p}  + (1-x) \ln \frac{1-x}{1-p} .$$
 We have the following large deviations estimate, for $s \in [0,1]$,
$$ \mathbb{P} \left( U_n \le s n \right) \le  \exp \{ - n \inf_{x \in [0, s]} H(x\,|\, p) \}.$$ 
\end{Lemma}

\begin{proposition}\label{propona}$\newline \vspace*{-0.5cm}$
 \begin{enumerate} 
\item[i)] Let $\nu$ be a vertex  with $|\nu|\ge 1$. The quantity
$$ \mathbb{P} \big(A(\nu) \;|\; h_{1}(\nu, {\rm fc}(\nu)) = x \big)$$
is a decreasing function of $x$, with $ x \ge 0$.
\item[ii)] $ \mathbb{P} \big(A(\nu) \;|\; h_{1}(\nu, {\rm fc}(\nu)) \le x \big) \ge \mathbb{P} \big(A(\nu)\big)$, for any $ x \ge 0$.
\end{enumerate}
\end{proposition}
\begin{proofsect}{Proof} Suppose $\{ {\rm fc}(\nu) = \overline{\nu}\}$. Given $\big\{ h_{1}(\nu,  \overline{\nu})=x\big\}$, the set of good vertices in $\Lambda_{ \overline{\nu}}$  is   a function of $x$. Denote this function by $\Tcal \colon \R^{+} \to \{$subset of vertices of $\Lambda_{ \overline{\nu}}\}$. A child of $ \overline{\nu}$, say $\nu_{1}$, is good if and only if 
$$h_{1}( \overline{\nu}, \nu_{1}) < \frac{h_{1}( \overline{\nu}, \nu)}{1 + x }.$$
 Hence the smaller $x$ is,  the more likely $\nu_{1}$ is good. This is true for any child of  $ \overline{\nu}$. As for descendants of $ \overline{\nu}$ at level  strictly greater than $|\nu|+2$, their status of being good is independent of $h_{1}( \nu, {\rm fc}(\nu))$. Hence $\Tcal(x) \supset \Tcal(y)$ for $x < y$.  This implies that the connected component of good vertices contining $ \overline{\nu}$ is larger if $\{h_{1}(\nu,  \overline{\nu})=x\}$ rather than $\{h_{1}(\nu,  \overline{\nu})=y\}$, for $x<y$.~~Hence
$$ 
\mathbb{P} \big(A(\nu) \;|\; h_{1}(\nu, {\rm fc}(\nu)) = x, \; {\rm fc}(\nu)= \overline{\nu}\big) \ge \mathbb{P} \big(A(\nu) \;|\; h_{1}(\nu, {\rm fc}(\nu)) = y,\;  {\rm fc}(\nu) =  \overline{\nu}\big), \qquad \mbox{for } x<y.
$$ 
Using symmetry we get i). In order to prove ii), use i) and the fact that the distribution of $h_{1}(\nu, \mbox{fc}(\nu))$ is stochastically larger that the conditional distribution of $h_{1}(\nu, \mbox{fc}(\nu))$ given $\{h_{1}(\nu, \mbox{fc}(\nu)) \le x\}$. \qed
\end{proofsect}
Denote by $[x]$ the largest integer smaller than $x$. 
\begin{theorem}\label{codaelle}
 For VRJP defined on $\Gcal_{b}$, with $b\ge 3$,  and  $s\in (0,1)$, we have
\begin{equation}\label{tailelle1}
 \mathbb{P} \bigl( l_{[s n]}  \ge n  \bigr) \le   \exp \Big\{ - [n/\zeta] \inf_{x \in [0, s]} H\Big(x \,\big{|}\, (1-\gamma_{b})\varphi_b \Big) \Big\}, 
\end{equation}
where $\gamma_{b}$ was defined in Lemma~\ref{return1}, and
\begin{equation}\label{fifi}
\varphi_{b} := \left( 1 - \e^{-b} \right) \left( 1 - \e^{-(b+1)} \right) \frac{b}{b+2}.
\end{equation}
\end{theorem}
\begin{proofsect}{Proof}  By virtue of Proposition~\ref{aind}   the sequence $\1_{A_{k \zeta}}$, with $ k \in \N$, consists  of i.i.d. random  variables. The random variable  $\sum_{j=1}^{[n/\zeta]} \1_{A_{j\zeta}}$ has binomial distribution   with parameters  $\big(\mathbb{P}\big(A(\rho)\big), [n/\zeta]\big)$.
 We define the event
 $$
 \begin{aligned}
 B_{j} := & \{\mbox{the first jump after $T_{j}$ is towards level $j+1$  and
 $L\big(X_{T_{j}}, T_{j+1}\big)<2$,}\\ 
 &\qquad \qquad \mbox{and $L\big(\mbox{par}(X_{T_{j}}), T_{j+1}\big)<2$}\}.
 \end{aligned}
 $$

Let $\Fcal_{t}$ be the smallest sigma-algebra defined by the collection $\{ X_{s},\, 0 \le s \le t \}$. For any stopping time $S$ define $\Fcal_{S} := \big\{ A \colon 
A \cap \{ S \le t \} \in \Fcal_{t} \big\}$. Now we show
\begin{equation}\label{typeb}
\begin{aligned}
 \mathbb{P} \left( B_{j} \, \mid\, \mathcal{F}_{T_{i -1 }} \right)  &\ge \left( 1 - \e^{-b} \right) \left( 1 - \e^{-(b+1)} \right) \frac{b}{b+2}=\varphi_b, 
\end{aligned} 
\end{equation}
\noindent where the inequality holds a.s.. In fact, by time $T_{i}$ the total weight of the parent of $X_{T_{i}}$ is stochastically smaller than $1+$ an exponential of parameter $b$, independent of $\Fcal_{T_{i-1}}$. Hence the probability that this total weight is less than 2 is larger than $1 - \e^{-b}$. Given this, the probability that    the first jump after $T_{i}$ is towards level $i+1$ is larger than $b/(b+2)$. Finally, the conditional probability that $T_{i+1}-T_{i} <1$ is larger than $1 - \e^{-(b+1)}$. This implies, together with $\zeta \ge 2$,  that the random variable  $\sum_{j=1}^{[n/\zeta]} \1_{B_{j}}$ is stochastically larger than a binomial$(n, \varphi_{b})$.   
For any  $i \in \N$,  and any vertex $\nu$ with $|\nu|=i\zeta$, set 
\begin{eqnarray*}
Z &:=& \min \Big( 1, \frac{h_{1}\big(\nu, \mbox{\rm par}(\nu)\big)}{1 + h_{1}\big( \mbox{\rm par}(\nu), \nu \big)}\Big)\\
E &:=& \{ X_{T_{i\zeta}} = \nu\} \cap\{ L(\mbox{par}(\nu), T_{i \zeta}) < 2\}. 
\end{eqnarray*}
We have 
$$ B_{i\zeta} \cap \{X_{T_{i\zeta}} = \nu\} = \{h_{1}(\nu,{\rm fc}(\nu)) <Z\} \cap E.$$
Moreover, the random variable $Z$  and the event $E$  are  both measurable with respect the sigma-algebra  
$$
 \widetilde{\Hcal}_{\nu} := \sigma\Big\{ P({\rm par}(\nu), \nu), \big\{  P(u,w) \colon \; u,w \notin \mbox{Vert}(\Lambda_{\nu}) \big\} \Big\}.
$$ 
Let $ f_Z$ be the density of $Z$ given $\{h_1(\nu, \rm{fc}(\nu)) < Z\}\cap E$. Using   \ref{propona}, ii), and the independence between $h_1(\nu, \rm{fc}(\nu))$ and $\widetilde{\Hcal}_{\nu}$,  we get
\begin{equation}\label{macv}
\begin{aligned}
\mathbb{P}\big(A_{i \zeta} \big{|} \,B_{i \zeta} \cap \{X_{T_{i \zeta}} = \nu\} \big)&= \mathbb{P}\big(A(\nu) \big{|} \,\{h_{1}(\nu, \mbox{fc}(\nu))< Z\}\cap E \big)\\
&=\int_0^\infty  \mathbb{P}\big(A(\nu) \big{|} \,\{h_{1}(\nu, \mbox{fc}(\nu))< z\}\big)f_Z(z) \d z \ge \mathbb{P}\big(A(\nu)\big)\\
&= \sum_{\nu \colon |\nu|= i \zeta} \mathbb{P}\big(A(\nu)\cap \{X_{T_{i \zeta}} = \nu\} \big)     =\mathbb{P}(A_{i \zeta}).
\end{aligned}
\end{equation}
 The first   equality in the last line of \eqref{macv}  is due  to symmetry.
Hence
\begin{equation}\label{posdip}
\mathbb{P}(A_{i \zeta}\big{|}B_{i \zeta} ) \ge \mathbb{P}(A_{i \zeta}).
\end{equation}
\noindent If $A_{k} \cap B_{k}$ holds then $k$ is a cut level. In fact, on this event, when the walk visits  level $k$ for the first time it jumps right away to  level $k+1$   and never visits level $k$ again. This happens  because  $X_{T_{k+1}} =$ fc$(X_{T_{k}})$ has a child which is  the root of an infinite subtree of good vertices. Moreover the total weights at $X_{T_{k}}$ and its parent are less than 2.  Define 
\[ 
e_n :=  \sum_{i=1}^{[n/\zeta]} \1_{A_{i\zeta}\cap B_{i \zeta}}.
 \]
By virtue of  \eqref{ret1}, \eqref{typeb}, \eqref{posdip} and Proposition~~\ref{aind} we have that $e_{n}$ is  stochastically larger than  a bin($[n/\zeta]$, $(1-\gamma_{b})\varphi_b$). Applying   Lemma~\ref{bind}, we have
\begin{equation*}
\begin{aligned}
\mathbb{P} \bigl( l_{[s n]} \ge n  \bigr) \le \mathbb{P} \bigl( e_{n} \le [s n]   \bigr) \le    \exp \Big\{ - [n/\zeta] \inf_{x \in [0, s]} H\Big(x \,\big{|}\, (1-\gamma_{b})\varphi_b \Big) \Big\}.\qquad \qquad \qquad \qquad \qed
\end{aligned}
\end{equation*}
\end{proofsect}
\noindent The function  $ H\big(x \,\big{|}\, (1-\gamma_{b})\varphi_b \big)$ is decreasing in the interval $(0, (1-\gamma_{b})\varphi_b)$. Hence for $ n > 1/\big((1-\gamma_{b})\varphi_b\big)$, we have $\inf_{x \in[0,1/n]} H\big(x \,\big{|}\, (1-\gamma_{b})\varphi_b \big) = H\big(1/n \,\big{|}\, (1-\gamma_{b})\varphi_b \big)$.
\begin{cor}\label{tail2ofl} For $n > 1/\big((1-\gamma_{b})\phi_{b}\big)$, by choosing $s = 1/n$ in Theorem~~\ref{codaelle},  we have
\begin{equation}\label{finaltailelle} 
\begin{aligned}
 \mathbb{P} \bigl( l_{1}  \ge n  \bigr) &\le   \exp \Big\{ - [n/\zeta]  \inf_{x \in [0, 1/n]}H\Big(x \,\big{|}\, (1-\gamma_{b})\varphi_b \Big) \Big\}\\
 &=   \exp \Big\{ - [n/\zeta]  H\Big(\frac 1n \,\big{|}\, (1-\gamma_{b})\varphi_b \Big) \Big\},
\end{aligned}
\end{equation}
where, from the definition of $H$ we have 
$$ \lim_{n \ti} H \Big(\frac 1n \,\big{|}\, (1-\gamma_{b})\varphi_b \Big)  = \ln \frac 1{1- (1-\gamma_{b})\varphi_b} >0.$$
\end{cor}


\section {  $\tau_1$ has finite $(2+\delta)$-moment}

The  goal of this section is to prove the finiteness of the $11/5$ moment of the first cut time. We adopt the following strategy
\begin{itemize}
\item first we prove the finiteness  of all moments for the number of vertices visited
by time $\tau_{1}$, then
\item we prove that  the total time spent at each of these sites has finite $12/5$-moment.
\end{itemize}
Fix $n\in \N$ and let
\begin{eqnarray*}
 &&\Pi_n :=  \mbox{ number of distinct vertices that $\mathbf{X}$ visits by time $T_n$}, \\
 &&\Pi_{n,k} := \mbox{number of distinct vertices that $\mathbf{X}$ visits at level $k$ by time $T_n$}.
\end{eqnarray*}

Let $T(\nu):= \inf \{ t \ge 0 \colon X_{t} = \nu\}$.  For any subtree $E$ of $\Gcal_{b}$, $ b \ge 1$,   define
\[ \delta(a, E) :=  
      \sup \left\{ t \colon \int_0^t \1_{\left\{ X_s \in E \right\}}\d s \le a \right\}. \] 
 The process  $X_{\delta(t,E)}$ is called the { \bfseries restriction } of  $\mathbf{X}$ to $E$.
 \begin{proposition}[{\bf Restriction principle (see \cite{DV2002})}] Consider VRJP $\mathbf{X}$ defined on a tree $\Jcal$ rooted at $\rho$. Assume this process is recurrent, i.e. visits each vertex infinitely  often, a.s.. Consider a subtree $\widetilde{\Jcal}$ rooted at $\nu$. Then the  process $X_{\delta(t,\widetilde{\Jcal})}$ is VRJP defined on $\widetilde{\Jcal}$. Moreover, for any subtree $\Jcal^{*}$ disjoint from $\widetilde{\Jcal}$, we have that 
 $X_{\delta(t,\widetilde{\Jcal})}$ and $X_{\delta(t,\Jcal^{*})}$ are independent. 
 \end{proposition}
\begin{proofsect}{Proof}
This principle follows directly from the Poisson construction and the memoryless property of  the exponential distribution.
\qed
\end{proofsect}

\begin{definition}
Recall that $P(x,y)$, with $x, y \in \mbox{Vert}\big(\Gcal_{b}\big) $ are the Poisson processes used to generate $\mathbf{X}$ on $\Gcal_{b}$. Let $\Jcal$ be a subtree of $\Gcal_{b}$. Consider VRJP $\mathbf{V}$ on $\Jcal$ which is generated by using $\big\{P(u,v)\colon u,v \in \mbox{Vert}(\Jcal)\big\}$, which is the same collection of Poisson processes used to generate the jumps of $\mathbf{X}$ from the vertices of $\Jcal$. We say that $\mathbf{V}$ is the {\bf extension}  of $\mathbf{X}$ in $\Jcal$. The processes $V_{t}$ and $X_{\delta(t,\Jcal)}$ coincide up to a random time, that is the total time spent by $\mathbf{X}$ in $\Jcal$.
\end{definition}
We construct an upper bound for  $\Pi_{n,k}$, with $ 2 \le k \le n-1$.
. Let $G(k)$ be the finite subtree of $\Gcal_b$ composed by all the vertices at level $i$ with $i \le k-1$, and the edges connecting them. Let $\mathbf{V}$ be the extension of $\mathbf{X}$ to $G(k)$. This process is recurrent, because is defined on a finite graph. The total number of first children at level $k-1$ is $b^{k-2}$, and we order them according to when they are visited by $\mathbf{V}$, as follows. Let $\eta_1$ be the first vertex at level $k-1$ to be visited by $\mathbf{V}$. Suppose we have defined $\eta_1, \ldots, \eta_{m-1}$. Let $\eta_m$ be the first child at level $k-1$ which does not belong to the set $\{\eta_1, \eta_2, \ldots, \eta_{m-1}\}$, to be visited. The vertices $\eta_i$, with $ 1 \le i \le b^{k-2}$ are determined by $\mathbf{V}$. All the other quantities and events such as $T(\nu)$ and $A(\nu)$, with $\nu$ running over the vertices of $\Gcal_b$, refer to the process $\mathbf{X}$. Define    
\[ f_{n}(k)  := 1 + b^2 \inf\{ m \ge 1\colon \1_{A(\rm{par}(\eta_m))} =1\}.   \]
Let $J:= \inf\{ n \colon T(\eta_n)= \infty\}$, if the infimum is over an empty set, let $J = \infty$. Suppose that $A(\eta_m)$ holds, then $\mathbf{X}$, after time $T(\eta_m)$, is forced to remain inside $\Lambda_{\eta_m}$, and never visits ${\rm fc}(\eta_m)$ again. This implies that $T(\eta_{m+1}) = \infty$. Hence, if $ J=m$ then $\bigcap_{i=1}^{m-1} (A({\rm par}(\eta_i)))^c$ holds, and $f_n(k) \ge 1+b^2(m-1)$. Similarly if $J= \infty$ then $f_n(k) = 	1+b^2 b^{k-2}= 1 + b^k$, which is an obvious upper bound for the number of vertices at level $k$ which are visited by $\mathbf{X}$. 
On the other  hand, if $J=m$ then the number of vertices at level $k$ which are visited by $\mathbf{X}$ is at most $ 1+(m-1)b^2$. In fact, the processes $\mathbf{X}$ and $\mathbf{V}$ coincide up to the random time when the former process leaves $G(k)$ and never returns to it. Hence if $T(\eta_i)<\infty$ then  $\mathbf{X}$ visited exactly $i-1$ distinct first children at level $k-1$ before time $T(\eta_i)$. On the event $\{J=m\}$ we have that $\{T(\eta_{m-1} < \infty\} \cap\{ T(\eta_m) = \infty\}$, hence exactly $m-1$ first children are visited at level $k-1$. This implies that at most $1+(m-1)b^2$ vertices at level $k$ are visited.

  We conclude that $f_{n}(k)$ overcounts the number of vertices at level $k$ which are visited, i.e.  $\Pi_{n,k} \le f_{n}(k)$. 

Recall that $h_{1}(\nu, {\rm fc}(\nu))$, being the minimum over a set of $b$ independent exponentials with rate 1, is distributed as an exponential with mean $1/b$.

\begin{Lemma}\label{cav1} For any  $m\in \N$, we have
\[ \mathbb{P} \bigl(\, f_{n}(k) \; > \; 1+ m b^2 \;  \bigr) \le  (\gamma_{b})^{m}. \]
\end{Lemma}
\begin{proofsect}{Proof}   Given $\bigcap_{i=1}^{m-1}(A({\rm par}(\eta_i))^c $ the distribution of $ h({\rm par}(\eta_m), \eta_m)$ is stochastically smaller than an exponential with mean $1/b$. Fix a set of vertices $\nu_i$ with $1 \le i \le m-1$ at level $k-1$ and each with a different parent. Given $\eta_i=\nu_i$ for $ i\le m-1$, consider the restriction of $\mathbf{V}$ to the finite subgraph obtained from $G(k)$ by removing each of the $\nu_i$ and ${\rm par}(\nu_i)$, with $ i\le m-1$. The restriction of $\mathbf{V}$ to this subgraph is VRJP, independent of $\bigcap_{i=1}^{m-1}(A({\rm par}(\eta_i))^c$, and the total time spent by this process in level $k-2$ is exponential with mean $1/b$. This total time is an upper bound for $ h({\rm par}(\eta_m), \eta_m)$. This conclusion is independent of our choice of the vertices $\nu_i$ with $1 \le i \le m-1$. Finally, using Proposition~\ref{propona} i), we have 
\begin{equation}
\begin{aligned}
\mathbb{P}\big(f_{n}(k)> 1 + mb^2\;|\; f_{n}(k)> 1 + (m-1)b^2 \big)&= \mathbb{P}\big( (A({\rm par}(\eta_{m})))^c \;|\; \bigcap_{i=1}^{m-1}(A({\rm par}(\eta_i))^c \big)\\
&\le \mathbb{P}\big( (A({\rm par}(\eta_{m})))^c \big) \le  \gamma_{b}.
\end{aligned}
\end{equation}  

\qed           \end{proofsect} 

Let $a_n, c_n$ be numerical sequences. We say that  $c_n= O(a_n)$ if $ c_n/ a_n $ is bounded.

\begin{Lemma}\label{finitep} For $p \ge 1$, we have $ \mathbb{E}\left [ \Pi^p_n  \right ]  = O(n^p). $ 
\end{Lemma}
\begin{proofsect}{Proof}   
Consider first the case  $p>1$. Notice that $\Pi_{n,0}=\Pi_{n,n}=1$. By virtue of Lemma~~\ref{cav1}, we have that 
$ \sup_{n} \mathbb{E}[f_{n}^{p}]<\infty.$
 By Jensen's inequality
\begin{equation}
\begin{aligned}
\mathbb{E} [ \Pi_n^p ] &= \mathbb{E} \left[ \left(2  +  \sum_{k=1}^{n-1} \Pi_{n,k}) \right)^p \right]  \le  n^p \; \mathbb{E}\left [ \sum_{k=1}^{n-1} \frac{\Pi_{n,k}^p}{n} + \frac{2^{p}}{n} \right] \le  n^p \; \mathbb{E}\left [ \sum_{k=1}^{n-1} \frac{f_{n}^p(k)}{n} + \frac{2^{p}}{n} \right]  \; = O(n^p). 
\end{aligned}
\end{equation}
\noindent   As for the case $p=1$,
\[ \mathbb{E} [ \Pi_n ] \le 2 + \sum_{k=1}^{n-1} \mathbb{E}[f_{n}(k)]= O(n). \]
\qed           \end{proofsect}
Let  
$$
\Pi   :=  \sum_{\nu}\1_{\large \{\nu \; \mbox{{\small is visited before time} }\,  \tau_{1}\}}.
$$
where the sum is over the vertices of $ \Gcal_b$. In words, $\Pi$ is the number of vertices visited before $\tau_1$.

\begin{Lemma}\label{finipi} For any $ p >0$ we have  $\mathbb{E} [ \Pi^p ] < \infty.$
\end{Lemma}
\begin{proofsect}{Proof}
By virtue of  Lemma~\ref{finitep}, $ \sqrt{\mathbb{E} \big[ \Pi_{n}^{2p} \big]}  \le C_{b,p}^{\ssup 1} n^p$, for some positive constant $C_{b,p}^{\ssup 1}$. Hence using Cauchy-Schwartz,
\begin{equation*}
\begin{aligned}
 \mathbb{E} \left[ \Pi^p \right] = \sum_{n= 1}^{\infty} \mathbb{E} \big [ \Pi_{n}^p \1_{ \left \{ l_1 = n \right \} }  \big ]&\le
\sum_{n= 1}^{\infty} \sqrt{\mathbb{E} \big [ \Pi_{n}^{2p} \big]\mathbb{P} (
l_1 \ge n )}\\  
&\le   C_{b,p}^{\ssup 1} \sum_{n= 1}^{\infty}  n^p   \exp \Big\{ - \frac 12 [n/\zeta]  H\Big(\frac 1n \,\big{|}\, (1-\gamma_{b})\varphi_b \Big) \Big\}    < \infty. 
\end{aligned}
\end{equation*}
\noindent In the last inequality we used Corollary~\ref{tail2ofl}. 
\qed           \end{proofsect}
 Next, we want to prove that the $12/5$-moment of $ L(\rho, \infty)$ is finite. We start with three intermediate results. The first two can be found in \cite{DV2004}. We include the proofs here for the sake of completeness.
\begin{Lemma}\label{btime}
Consider   VRJP on $\{ 0, 1 \}$, which starts at $1$, and with initial weights $ a_0 = c$ and $ a_1 =1$.  Define 
\[
 \xi(t):= \inf \Bigl\{ s \colon L(1,s)= t \Bigr\}. 
 \] 
We have
 \begin{equation}\label{boundtimesq}
 \sup_{ t \ge 1} \, \mathbb{E} \left[ \left( \frac{L(0, \xi(t) )}{t} \right)^3 \right] = c^{3}+ 3c^2+3c. 
 \end{equation}
\end{Lemma}
\begin{proofsect}{Proof}
 We have $L(0, \xi(t + \d t)) = L(0, \xi(t)) + \chi \eta$, where $ \chi$ is a Bernoulli which takes value 1 with probability $ L(0, \xi(t))\d t$,  and $ \eta$ is exponential with mean $1/t$. Given $L\big(0, \xi(t)\big)$, the random variables  $ \chi$ and $ \eta$ are independent. 
 Hence
\[ \mathbb{E} \Bigl[ L(0, \xi(t+\d t) )  \Bigr] -  \mathbb{E} \Bigl[ L(0, \xi(t) ) \Bigr ] =
\frac{\mathbb{E} [ L(0, \xi(t) )  ]}{t}\d t, \] 
i.e. $\mathbb{E} [L(0, \xi(t) )]$ is solution of the equation $y^{'}(t) = y(t)/t$, with initial condition $ y(1)= c$ (see \cite{DV2002}). Hence 
\[   \mathbb{E} [ L(0, \xi(t) )  ] = ct. \] 
Similarly
\begin{equation*}
\begin{aligned}
&  \mathbb{E} \Bigl[ L(0, \xi(t+\d t) )^2 \Bigr ]  \\
& \quad = \mathbb{E} \Bigl[ L(0, \xi(t) )^2 \Bigr]  +  2 \mathbb{E} \Bigl[ L(0, \xi(t) )  \mathbb{E} \Bigl[ \chi \; \mid \; L(0, \xi(t) ) \Bigr ] \Big] \mathbb{E} [ \eta ]  + \; \mathbb{E} \Bigl[ \chi^2 \; \mid \; L(0, \xi(t) ) \Bigr ] \mathbb{E} [ \eta^2 ]  \\
 &\quad = \mathbb{E} \Bigl[ L(0, \xi(t) )^2 \Bigr ] + (2/t) \mathbb{E} \Bigl[ L(0, \xi(t) )^2 \Bigr]\d t + (2/t^2)\mathbb{E} \Bigl[ L(0, \xi(t) ) \Bigr]\d t \\
 &\quad = \mathbb{E} \Bigl[ L(0, \xi(t) )^2 \Bigr ] + (2/t) \mathbb{E} \Bigl[ L(0, \xi(t) )^2 \Bigr]\d t + (2c/t)\d t. \\
\end{aligned}
\end{equation*}
\noindent Thus $ \mathbb{E} \Bigl[ L(0, \xi(t) )^2 \Bigr ]$ satisfies the equation $y^{'} = (2/t) y + (2c/t) $, with $y(1)=c^2$. Then,
\[ \mathbb{E} \Bigl[ L(0, \xi(t) )^2 \Bigr ] = - c + \left( c^2 + c \right)  t^2. \]
Finally, reasoning in a similar way, we get that $ \mathbb{E} \Bigl[ L(0, \xi(t) )^3 \Bigr ]$ satisfies the equation $y^{'} = (3/t) y + 6(c^{2}+c) $, with $y(1)=c^3$. Hence,
\[ \mathbb{E} \Bigl[ L(0, \xi(t) )^3 \Bigr ] = - 3(c^{2}+c)t  + \left( c^{3}+ 3c^2 +3 c \right)  t^3. \]
 Divide both sides by $t^3$, and use the fact that  $c>0$ to  get \eqref{boundtimesq}.
\qed           \end{proofsect}



%
A ray $\sigma$ is a subtree of $\Gcal_{b}$ containing exactly one vertex of each level of $\Gcal_{b}$.
  Label the vertices of this ray using $\{\sigma_{i},\,  i \ge 0\}$,  where  $\sigma_{i}$ is  the unique vertex at level $i$ which belongs to $\sigma$. Denote by $\Sym$ the collection of all rays of $\Gcal_{b}$.

\begin{Lemma}\label{spbe}
For any ray $\sigma$, consider VRJP $\; \mathbf{X}^{\ssup \sigma}:= \{ X^{\ssup \sigma}_{t}, \, t \ge 0\}$, which is the extension of $\mathbf{X}$ to $\sigma$. Define
 \begin{eqnarray*}
 T^{\ssup \sigma}_{n} &:=& \inf \{ t>0 \colon X^{\ssup \sigma}_{t} = \sigma_{n} \},\\
L^{\ssup \sigma} ( \sigma_{i},t )  &:=& 1 + \int_{0}^{t}\1_{\{X_{s}^{\ssup{\sigma}}= \sigma_{i}\}} \d s.
 \end{eqnarray*}
  We have that 
\begin{equation}\label{upbsq}
 \mathbb{E} \big[ L^{\ssup \sigma}(\sigma_{0}, T^{\ssup \sigma}_n)^3 \big] \le (37)^n. 
\end{equation}
\end{Lemma}
\begin{proofsect}{Proof}  By the tower property of conditional expectation,
\begin{equation}\label{reiterate1}
 \mathbb{E} \left[ \bigl(L^{\ssup \sigma} ( \sigma_{0},T^{\ssup \sigma}_n ) \bigr)^3 \right] = 
\mathbb{E} \left[ \bigl(L^{\ssup \sigma} ( \sigma_{1}, T^{\ssup \sigma}_n )\bigr)^3  \mathbb{E} \left[  \left (\frac{L^{\ssup \sigma} \left( \sigma_{0}, T^{\ssup \sigma}_n \right)}{L^{\ssup \sigma} \left( \sigma_{1},T^{\ssup \sigma}_n\right)} \right)^3 \Big{|} \, L^{\ssup \sigma} \left( \sigma_{1},T^{\ssup \sigma}_n \right) \right]\right]. 
\end{equation}
At this point  we focus on  the process restricted to $\{ 0, 1 \}$. 
 This restricted process is  VRJP which starts at $1$, with initial weights $a_{1} = 1 $, and $a_{0} =   1 + h_{1}(\sigma_{0}, \sigma_{1})$ and $\sigma_{0} = \rho$.   By applying Lemma~\ref{btime}, and using the fact that $ h_{1}(\sigma_{0}, \sigma_{1})$ is exponential with mean 1, we have
\begin{equation}\label{upt1}
\begin{aligned}
\mathbb{E} \left[  \left(\frac{L^{\ssup \sigma} \left( \sigma_{0}, T^{\ssup \sigma}_n \right)}{L^{\ssup \sigma} \left( \sigma_{1},T^{\ssup \sigma}_n \right)} \right)^3 \Big{|}  L^{\ssup \sigma} \left( \sigma_{1},T^{\ssup \sigma}_n \right) \right] &\le \mathbb{E}\big[3(1+h_{1}(\sigma_{0}, \sigma_{1})) + (1+h_{1}(\sigma_{0}, \sigma_{1}))^{2} + (1+h_{1}(\sigma_{0}, \sigma_{1}))^{3} \big]\\ 
&= 37.
\end{aligned}
\end{equation}
\noindent Then
\begin{equation}
\begin{aligned}
\mathbb{E} \left[ \bigl( L ( \sigma_{0},T_n ) \bigr)^3 \right] &=\mathbb{E} \left[  \mathbb{E} \left[  \left(\frac{L^{\ssup \sigma} \left( \sigma_{0}, T^{\ssup \sigma}_n \right)}{L^{\ssup \sigma} \left( \sigma_{1},T^{\ssup \sigma}_n \right)} \right)^3 \Big{|}  L^{\ssup \sigma} \left( \sigma_{1},T^{\ssup \sigma}_n \right) \right] \bigl( L^{\ssup \sigma} ( \sigma_{1},T^{\ssup \sigma}_n ) \bigr)^3 \right]\\
&\le 37\, \mathbb{E} \left[ \bigl(L^{\ssup \sigma} ( \sigma_{1}, T^{\ssup \sigma}_n )\bigr)^3 \right]. 
\end{aligned}
\end{equation}
The Lemma  follows by recursion and restriction principle.   
\qed           \end{proofsect}
Next, we prove that
\begin{equation}\label{spchlac}
 L(\rho, T(\sigma_{n}))  \le L^{\ssup \sigma}(\sigma_{0}, T^{\ssup \sigma}_{n}).
\end{equation}
In fact, we have equality if $T(\sigma_{n})<\infty$, because the restriction and the extension of $\mathbf{X}$ to $\sigma$  coincide during the time interval $[0, T(\sigma_{n})]$. If $T(\sigma_{n})=\infty$, it means that $\mathbf{X}$ left the ray $\sigma$ at a  time $s<T^{\ssup \sigma}_{n}$. Hence 
$$  L(\rho, T(\sigma_{n})) = L^{\ssup \sigma}(\sigma_{0}, s) \le L^{\ssup \sigma}(\sigma_{0}, T^{\ssup \sigma}_{n}).$$
Hence, for any $\nu$, with $|\nu|=n$, we have
\begin{equation}\label{timerho}
 \mathbb{E} \big[ L(\rho, T(\nu))^3 \big] \le (37)^n.
\end{equation}
\begin{Lemma}\label{uptinf} $ \mathbb{E} \left [ \left( L (\rho , \infty ) \right)^{12/5}  \right ] < \infty. $
\end{Lemma}
\begin{proofsect}{Proof}   Recall the definition of $A(\nu)$ from \eqref{sub1} and set
$$ D_{k} := \bigcup_{\nu \colon |\nu|=k-2} A(\nu).$$
 If $A(\nu)$ holds, after the first time the process hits the first child of $\nu$, if this ever happens,  it will never visit $\nu$   again, and will not increase the local time spent at the root. Roughly,~~our strategy~~is to use the extensions on paths to give an upper bound of the total time spent at the root by
 time $T_{k}$ and show that the probability that $\bigcap_{i=1}^k D^{c}_{i}$ decreases quite fast in $k$.

   Using the independence between disjoint collections of Poisson processes, we infer that $A(\nu)$, with $|\nu|=k-2$ are independent. In fact each $A(\nu)$ is determined by the Poisson processes attached to pairs of vertices in $\Lambda_{\nu}$. Hence
   \begin{equation}\label{pippo2}
 \begin{aligned}
\mathbb{P}(D^{c}_{k} ) \le (\gamma_{b})^{b^{k-2}}
 \end{aligned}
 \end{equation}
 Define $d = \inf\{ n \ge 1 \colon \1_{D_{n}} =1\}$. 
Fix $ k \in \N$.  On the set $\{d=k\}$, define  $\overline{\mu}$ to be one~of~the~first children at level $k-1$ such that $A($par$(\mu))$ holds.   On $\{T(\overline{\mu}) < \infty\} \cap \{d=k\}$, we clearly have $L(\rho, \infty)= L(\rho, T(\overline{\mu}))$. On the other hand, on $\{T(\overline{\mu}) < \infty\} \cap \{d=k\}$, we have that, after the process reaches  $\overline{\mu}$ it will never return to the root. Hence 
$$ L(\rho, \infty) = 1 + \int_{0}^{T(\overline{\mu})} \1_{\{X_{u} = \rho\}} \d u + \int_{T(\overline{\mu})}^{\infty} \1_{\{X_{u} = \rho\}} \d u = 1 + \int_{0}^{T(\overline{\mu})} \1_{\{X_{u} = \rho\}} \d u = L(\rho, T(\overline{\mu})).$$
Using this fact, combined with
$$  L(\rho, T(\overline{\mu})) \le \sum_{\nu \colon |\nu| = k-2}L(\rho, T({\rm fc}(\nu)),$$
and   $\1_{\{d=k\}} \le \1_{\{d > k-1\}} \le \1_{D^{c}_{k-1}}$, we have 
\begin{equation}\label{pippo1}
\begin{aligned}
L(\rho, \infty)\1_{\{d=k\}} = L(\rho, T(\overline{\mu})) \1_{\{d=k\}}  &\le \Big(\sum_{\nu \colon |\nu| = k-2}L\big(\rho, T({\rm fc}(\nu))\big)\Big)\1_{\{d=k\}}\\ 
&\le \Big(\sum_{\nu \colon |\nu| = k-2}L\big(\rho, T({\rm fc}(\nu))\big)\Big)\1_{D^{c}_{k-1}}. 
\end{aligned}
\end{equation}
Using \eqref{pippo1}, Holder's inequality (with $p= 5/4$) and \eqref{pippo2} we have
\begin{equation*}
\begin{aligned}
&\mathbb{E}  \left [ \left( L (\rho , \infty ) \right)^{12/5}  \right ] = \sum_{k=1}^{\infty} \mathbb{E} \left [ \big( L (\rho , \infty ) \big)^{12/5} \1_{\{d=k\}} \right]= \sum_{k=1}^{\infty} \mathbb{E} \left [ \big( L (\rho , \infty )  \1_{\{d=k\}} \big)^{12/5} \right]\\
&\le  \sum_{k=1}^{\infty} \mathbb{E}  \left [ \left( \sum_{\nu \colon |\nu| = k-2}L\big(\rho, T({\rm fc}(\nu))\big)\1_{D^{c}_{k-1}} \right)^{12/5}\right]
\le    \sum_{k=1}^{\infty} \mathbb{E} \left [
\left ( \sum_{\nu \colon |\nu| = k-2}L\big(\rho, T({\rm fc}(\nu))\big)\right)^{3} \right ]^{4/5} ( \gamma_{b})^{ b^{k-3}/5}\\  
&\le    \sum_{k=1}^{\infty} \mathbb{E} \left [ 
\left (  \sum_{\nu \colon |\nu| = k-2}L\big(\rho, T({\rm fc}(\nu))\big)\right)^{3} \right ] ( \gamma_{b} )^{ b^{k-3}/5} \qquad \mbox{(using $L(\rho, t) \ge 1$)} \\ 
&\le    \sum_{k=1}^{\infty}  \left (  b^{2k}
\sum_{\nu \colon |\nu| = k-2} \mathbb{E}\big[L\big(\rho, T({\rm fc}(\nu))\big)^{3}\big] \right ) ( \gamma_{b} )^{b^{k-3}/5}  \qquad \mbox{(by Jensen)}\\
&\le   \sum_{k=1}^{\infty}b^{3k} (37)^{ k} (\gamma_{b})^{b^{k-3}/5} < \infty. \\
\end{aligned}
\end{equation*} 
\qed           \end{proofsect}

\begin{Lemma}\label{randbt}
For $\nu \neq \rho$, there exists a random variable $\Delta_{\nu}$ which is $\sigma\big\{P(u,v)\colon u,v \in \mbox{Vert}(\Lambda_{\nu})\big\}$-measurable, such that
\begin{itemize}
\item[i)] $L(\nu, \infty) \le \Delta_{\nu}$, and
\item[ii)] $\Delta_{\nu}$ and $ L(\rho, \infty)$ are identically distributed.  
\end{itemize}
\end{Lemma}
\begin{proofsect}{Proof}
       Let $ \widetilde{\mathbf{X}}:= \{ \widetilde{X}_{t},\, t \ge 0\}$ be the extension of $ \mathbf{X} $ on $\Lambda_{\nu}$. Define 
\[ \Delta_{\nu} := 1 + \int_{0}^{\infty} \1_{\{ \widetilde{X}_t = \nu \}}\d t. \]
By  construction, this random variable satisfies i) and ii) and is $\sigma\big\{P(u,v)\colon u,v \in \mbox{Vert}(\Lambda_{\nu})\big\}$-measurable. 
\qed           \end{proofsect}

\begin{theorem}\label{boundfirstc}  $ \mathbb{E} \left [ (\tau_1)^{11/5} \right ]  < \infty. $
\end{theorem}
\begin{proofsect}{Proof} Suppose we relabel the vertices that have been  visited by time $\; \tau_1,$ using $ \theta_1, \theta_2, \ldots, \theta_{\Pi}$, where vertex $\nu$ is labeled $\theta_k$ if there are exactly $ k-1$ distinct vertices that have been visited before $ \nu$. Notice that $ \Delta_{\nu}$ and $ \{ \theta_k = \nu \} $ are independent,
because they are determined by disjoint non-random sets of Poisson processes ($ \Delta_{\nu}$ is  $\sigma\big\{P(u,v)\colon u,v \in \mbox{Vert}(\Lambda_{\nu})\big\}$-measurable).   As the variables $\Delta_{\nu}$, with $\nu \in$ Vert($\Gcal_{b})$, share the same distribution, for any $ p > 0$,  we have 
\[ \mathbb{E} [ \Delta_{\theta_k}^p ] = \mathbb{E} [ \Delta_{\nu}^p ] = \mathbb{E} [ L(\rho, \infty )^p]. \]
  By Jensen's and   Holder's (with $p = 12/11$) inequalities,   Lemma~\ref{randbt}  i) and ii), and Lemma~\ref{uptinf}, we have 
\begin{equation*}
\begin{aligned}
 \mathbb{E }  \left [ (\tau_1)^{11/5} \right ] &\le \mathbb{E } \left [ \left(\sum_{k=1}^{\Pi }
\Delta_{\theta_k} \right)^{11/5} \right ] \le \mathbb{E } \left [ \Pi^{(11/5) -1}\sum_{k=1}^{\Pi }
(\Delta_{\theta_k})^{11/5}  \right ]\\ &= \mathbb{E } \left [ \sum_{k=1}^{\infty } 
\Delta_{\theta_k}^{11/5}  \; \Pi^{6/5}\; \1_{\left \{\Pi \ge k  \right \} } \right ]  \le \sum_{k=1}^{\infty } \mathbb{E} \left [ \Delta_{\theta_k}^{12/5} \right ]^{11/12} \mathbb{E } \left [ \Pi^{72/5}\, \1_{\left ( \Pi \ge k  \right ) } \right ]  
^{1/12}  \\  
&\le C_b^{\ssup 3} \sum_{k=1}^{\infty }  \mathbb{E } \left [ \Pi^{144/5} \right]^{1/24}  \mathbb{P}\left ( \Pi \ge k  \right )^{1/24}  \qquad \mbox{(by Cauchy-Schwartz and Lemma~\ref{uptinf})} \\ 
&\le  C_b^{\ssup 4}  \sum_{k=1}^{\infty } \mathbb{P}  \left ( \Pi \ge k  \right )^{1/24}, \qquad \qquad \qquad \mbox{(by Lemma~\ref{finipi})},
\end{aligned}
\end{equation*}
for some positive constants $C_b^{\ssup 3}$ and $C_b^{\ssup 4}$.
It remains to prove the finiteness of the last sum. We use the fact
\begin{equation}\label{tailpi}
\lim_{k \to \infty} k^{48} \mathbb{P}( \Pi \ge k) = 0.
\end{equation}
The previous limit is a consequence of the well-known formula
\begin{equation}
\sum_{k=1}^{\infty} k^{48} \mathbb{P}( \Pi \ge k) = \mathbb{E}[ \Pi^{49}],
\end{equation}
and the finiteness of $ \mathbb{E}[ \Pi^{49}]$ by virtue of Lemma~\ref{finipi}. $$
\sum_{k=1}^{\infty } \mathbb{P}  \left ( \Pi \ge k  \right )^{1/24} = \sum_{k=1}^{\infty} \frac 1{k^{2}}\Big( k^{48} \mathbb{P}(\Pi \ge k)\Big)^{1/24}< \infty.
$$
\qed
\end{proofsect}

\begin{Lemma}\label{primloc}
 $\sup_{x \in [1,2]}\mathbb{E}[(L(\rho,\infty))^{12/5}\;|\; L(\rho,T_{1}) =x] < \infty$. 
\end{Lemma}
\begin{proofsect}{Proof} Using \ref{randbt}, and the fact that $\Delta_{X_{T_{1}}}$ is independent of $L(\rho, T_{1})$, we have 
\begin{equation}
\begin{aligned}
\sup_{x \in [1,2]}\mathbb{E}[(L(X_{T_{1}},\infty))^{12/5}\;|\; L(\rho,T_{1}) =x] &\le  \sup_{x \in [1,2]}\mathbb{E}[(\Delta_{X_{T_{1}}})^{12/5}\;|\; L(\rho,T_{1}) =x]\\
&= \mathbb{E}[(\Delta_{X_{T_{1}}})^{12/5}] = \mathbb{E}[(L(\rho,\infty))^{12/5}]<\infty.
\end{aligned}
\end{equation}

Given $L(\rho,T_{1}) =x$, the process $\mathbf{X}$ restricted to $\{\rho, X_{T_{1}}\}$ is VRJP which starts from $ X_{T_{1}}$, with initial weights $a_{\rho} =x$ and $1$  on $ X_{T_{1}}$.  This process runs up to the last visit of $\mathbf{X}$ to one of these two vertices. Using Lyapunov inequality, i.e. $\mathbb{E}[Z^{q}]^{1/q} \le \mathbb{E}[Z^{p}]^{1/p}$  whenever $ 0<q \le p$,     Lemma~\ref{spbe}, and the fact $x \ge 1$,  we have 
\begin{equation}
\begin{aligned}
 \mathbb{E}\Big[&\Big(\frac{L(\rho,T_{n})}{L(X_{T_{1}},T_{n})}\Big)^{12/5}\;|\; L(X_{T_{1}},T_{n}), \{L(\rho,T_{1})=x\}\Big]\\  
 &\le   \mathbb{E}\Big[\Big(\frac{L(\rho,T_{n})}{L(X_{T_{1}},T_{n})}\Big)^{3} \;|\; L(X_{T_{1}},T_{n}), \{L(\rho,T_{1})=x\}\Big]^{4/5}\\
& \le   (x^{3} + 3 x^{2} + 3x)^{4/5}  \le x^{3} + 3 x^{2} + 3x.
\end{aligned}
\end{equation}
Finally
\begin{equation}
\begin{aligned}
\mathbb{E}[(L(\rho,T_{n}))^{12/5}\;|\; L(\rho,T_{1}) =x] &= \mathbb{E}\Big[\Big(\frac{L(\rho,T_{n})}{L(X_{T_{1}},T_{n})}\Big)^{12/5}(L(X_{T_{1}},T_{n}))^{12/5} \;|\; L(\rho,T_{1}) =x\Big]\\
&\le (x^{3} + 3 x^{2} + 3x)  \mathbb{E}\Big[(L(X_{T_{1}},T_{n}))^{12/5} \;|\; L(\rho,T_{1}) =x\Big]\\
&\le (x^{3} + 3 x^{2} + 3x)  \mathbb{E}\Big[(L(X_{T_{1}}, \infty))^{12/5} \;|\; L(\rho,T_{1}) =x\Big]\\
&\le (x^{3} + 3 x^{2} + 3x)  \mathbb{E}[(L(\rho,\infty))^{12/5}].
\end{aligned}
\end{equation}
By sending $n \ti$ and taking the suprema over $x \in [1,2]$ we get 
$$ 
\sup_{x \in [1,2]}\mathbb{E}[(L(\rho,T_{n}))^{12/5}\;|\; L(\rho,T_{1}) =x] \le 26 \mathbb{E}[(L(\rho,\infty))^{12/5}]<\infty.
$$
\qed
\end{proofsect}
\begin{theorem} $ \sup_{x \in [1,2]} \mathbb{E }  \left [ (\tau_1)^{11/5} \; \mid L\bigl( \rho, T_1 \bigr)=x \right ] < \infty.$
\end{theorem}
\begin{proofsect}{Proof}         
Label the vertices at level 1 by  $\mu_1, \mu_2, \ldots, \mu_b$.
 Let $\tau_1(\mu_i)$ be 
the first cut time of the extension of $\mathbf{X}$ on $\Lambda_{\mu_i}$.  This extension is VRJP on $ \Lambda_{\mu_i}$  with initial weights 1, hence we can apply  Theorem~\ref{boundfirstc} to get   
\begin{equation}
 \mathbb{E} [\left( \tau_1(\mu_i) \right)^{11/5} ] < \infty.
 \end{equation}
  Hence, it remains to prove that  for $x \in [1,2]$ 
\[
\begin{aligned}
 &\mathbb{E} \left[ \big( \tau_1\big)^{11/5}  \,|\,  L\bigl( \rho, T_1 \bigr)=x \right] \le  \mathbb{E}\Big[  \Big(L(\rho, \infty)+\max_i \tau_{1}(\mu_{i})\Big)^{11/5} \,\Big{|}\,  L\bigl( \rho, T_1 \bigr)=x \Big]\\ &\le  \mathbb{E}\Big[  \Big(L(\rho, \infty)+\sum_{i=1}^{b}\tau_{1}(\mu_{i})\Big)^{11/5} \,\Big{|}\,  L\bigl( \rho, T_1 \bigr)=x  \Big]\\
&\le  (b+1)^{11/5 -1} \mathbb{E}\Big[  \Big(L(\rho, \infty)\Big)^{11/5}  \,|\,  L\bigl( \rho, T_1 \bigr)=x \Big] +(b+1)^{11/5}  \mathbb{E}\Big[ \big((\tau_{1}(\mu_{1})\big)^{11/5}\Big]< \infty,
 \end{aligned}
 \]
 where we used Jensen's inequality, the independence of $\tau(\mu_{i})$ and $T_{1}$ and Lemma~\ref{primloc}. In fact, as $ L\bigl( \rho, \infty \bigr) \ge 1$ , we have 
 $$ \mathbb{E}[(L(\rho,\infty))^{11/5}\;|\; L(\rho,T_{1}) =x]  \le \mathbb{E}[(L(\rho,\infty))^{12/5}\;|\; L(\rho,T_{1}) =x] < \infty.$$ 
\qed           \end{proofsect}

\section{Splitting the path into  one-dependent pieces} 
\noindent  Define $Z_i = L (X_{\tau_i}, \infty)$, with $ i \ge 1 $. 
\begin{Lemma}\label{process}
The process $Z_{i}$, with $i \ge 1$ is a homogenous Markov chain with state space $[1,2]$.
\end{Lemma}
\begin{proofsect}{Proof}
Fix $n \ge 1$.  On $\{Z_{n}=x\} \cap\{ X_{\tau_{n}} = \nu\}$ the random variable $Z_{n+1}$ is  determined by the variables $\{ P(u,v), u,v \in \Lambda_{\nu}, u \neq \nu\}$. In fact these  Poisson processes, on the set $\{Z_{n}=x\} \cap\{ X_{\tau_{n}} = \nu\}$,  are the only ones used to generate  the jumps of the process $\{X_{T(\rm{fc}(\nu)+t}\}_{t \ge 0}$.  Let $E_{1}, E_{2}, \ldots, E_{n-1}, E_{n+1}$ be Borel subsets of $[0,1]$. Conditionally on $\{Z_{n}=x\} \cap\{ X_{\tau_{n}} = \nu\}$, the two events $\{Z_{n+1} \in E_{n+1}\}$ and $\{Z_{1} \in E_{1}, Z_{2} \in E_{2}, \ldots Z_{n-1} \in E_{n-1}\}$ are independent because are determined by disjoint collections of Poisson processes. By symmetry 
$$ \mathbb{P}(Z_{n+1}\in E_{n+1}\;|\; \{Z_{n}=x\} \cap\{ X_{\tau_{n}} = \nu\})$$
does not depend on $\nu$. Hence
$$ 
\begin{aligned}
&\mathbb{P}(Z_{n+1}\in E_{n+1}\;| \; Z_{1} \in E_{1}, Z_{2} \in E_{2}, \ldots Z_{n-1} \in E_{n-1}, Z_{n}=x)\\
&= \sum_{\nu}  \mathbb{P}(Z_{n+1}\in E_{n+1}\;| \; Z_{1} \in E_{1},  \ldots, Z_{n-1} \in E_{n-1}, Z_{n}=x,  X_{\tau_{n}} = \nu ) \mathbb{P}( X_{\tau_{n}} = \nu \;|\;  Z_{1} \in E_{1},  \ldots,  Z_{n}=x)\\
&= \mathbb{P}(Z_{n+1}\in E_{n+1}\;| \;  Z_{n}=x,  X_{\tau_{n}} = \nu ) = \mathbb{P}(Z_{n+1}\in E_{n+1}\;| \;  Z_{n}=x).
\end{aligned}
$$
 This implies that $\mathbf{Z}$ is a Markov chain. The self-similarity property of $\Gcal_{b}$ and $\mathbf{X}$ yields the homegeneity.
\qed
\end{proofsect}

From the previous proof, we can infer that given $Z_i = x $, the random vectors 
$(\tau_{i+1}- \tau_i,\, l_{i+1} - l_i)$ and $ (\tau_{i}- \tau_{i-1},\, l_i - l_{i-1})$, are  independent.

\begin{proposition} 
\begin{eqnarray}
&&\label{supera1}\sup_{i \in \N} \sup_{x \in [1,2]} \mathbb{E} \Big[ \big(\tau_{i+1} - \tau_{i}\big)^{11/5} \, \big{|}\, Z_{i} =x \Big] < \infty\\
&&\label{supera2}\sup_{i \in \N} \sup_{x \in [1,2] } \mathbb{E} \left[ \bigl(l_{i+1} - l_i \bigr)^{11/5} \; \mid \; Z_i=x \right] < \infty. 
\end{eqnarray}
\end{proposition}
\begin{proofsect}{Proof}
We only prove \eqref{supera1}, the proof of  \eqref{supera2} being similar. Define $C:= \big\{ X_{t} \neq \rho, \, \forall t > T_{1} \big\}$ and  fix a vertex $\nu$.  Notice that by the self-similarity property of $\Gcal_{b}$, we have
\begin{equation*}
 \mathbb{E} \Big[  ( \tau_{i+1} - \tau_i)^{11/5}  \; \mid \; \{ Z_i=x \} \cap \{  X_{\tau_i} = \nu \} \Bigr] = \mathbb{E} \Big[ (\tau_{1})^{11/5} | \{ L(\rho, T_{1}) = x \} \cap C \Big].
\end{equation*}
By the proof of Lemma~\ref{return},  we have that 
\begin{equation}\label{return3}
\inf_{1 \le x \le 2} \mathbb{P} \big( C \, \big{|} \, L(\rho, T_{1}) = x \big) \ge \frac b{b+x} \mathbb{P}(A_{1}) \ge (1-\gamma_{b})\frac b{b+2}>0.
\end{equation}
Hence
\begin{equation*}
\begin{aligned}
\sup_{x \colon x \in [1,2]} \mathbb{E} & \Big[ (\tau_{1})^{11/5} \,\big{|} L(\rho, T_{1}) =x \Big]\\
&\ge \sup_{x \colon x \in [1,2]} \mathbb{E}  \Big[ (\tau_{1})^{11/5} \,\big{|}\{ L(\rho, T_{1}) =x\} \cap C \Big] \mathbb{P}(C\, |\, L(\rho, T_{1}) =x)\\
&\ge (1-\gamma_b) \frac b{b+2} \sup_{x \colon x \in [1,2]} \mathbb{E}  \Big[ (\tau_{1})^{11/5} \,\big{|}\{ L(\rho, T_{1}) =x\} \cap C \Big]\\
&\ge (1- \gamma_b)\frac b{b+2}  \sup_{x \colon x \in [1,2]} \mathbb{E} \Big[ (\tau_{i+1} - \tau_{i})^{11/5} \, | \, \{ Z_{i} = x\} \cap \{ X_{\tau_{i}} = \nu\} \Big]
\end{aligned}
\end{equation*}
Hence
\begin{equation*}
\begin{aligned}  
&\mathbb{E}\Bigl[  (\tau_{i+1} - \tau_i)^{11/5} \mid \{ Z_i=x \} \cap \{  X_{\tau_i} = \nu \}
\Bigr] \le   \frac {b+2}{b(1- \gamma_{b})}\sup_{1\le x \le 2 } \mathbb{E} \Bigl[ (\tau_{1})^{11/5} \mid \{ L(\rho,T_1  = x \}  \Bigr]. 
\end{aligned}
\end{equation*}
\qed           \end{proofsect}

Next we prove that $\mathbf{Z} $ satisfies the Doeblin condition.

\begin{Lemma}
There exists a  probability measure $\phi(\cdot)$  and  $ 0< \lambda \le  1 $, such that for every Borel subset $B$ of $[1,2]$, we have 
\begin{equation}\label{homzeta}
 \mathbb{P} \big( Z_{i+1} \in B \; \mid Z_i = z \big) \; \ge \;  \lambda \; \phi(B)\qquad  \forall \; z \in [1,2]. 
\end{equation}
\end{Lemma}
\begin{proofsect}{Proof}
As $Z_{i}$ is homogeneous, it is enough to prove \eqref{homzeta} for $i=1$. In this proof we show that the distribution of
$Z_{2}$ is absolutely continuous and we compare it to $1+$ an exponential  with parameter 1 conditionated on being less than 1. The analysis is technical because $Z_{i}$ depend on the behaviour of the whole process $\mathbf{X}$.
Our goal is to find a lower bound for 
\begin{equation}\label{crib}
 \mathbb{P} \Big( Z_2 \; \in \; (x,y) \; \big{|} \;  Z_1 =z   \Big), \qquad \mbox{ with $z \in [1,2]$}.
\end{equation}
Moreover, we require that this lower bound is independent of $z \in [1,2]$.

Fix $\eps\in (0,1)$.  Our first goal is to find a lower bound for the probability of the event $\{ Z_2  \in  (x,y), \;  Z_1 \in I_{\eps}(z) \}$, where $I_{\eps}(z):= (z- \eps, z+\eps)$. 
Fix $z \in [1,2]$  and consider the function
\begin{equation}\label{thefun}
   \e^{-(b+u)(t-1) } - (b+1)\e^{-(b+2)} \e^{-(t-1)}. 
\end{equation}
Its derivative with respect $t$ is 
$$ (b+1)\e^{-(b+2)- (t-1)}- (b+u)\e^{-(b+u)(t-1)},$$
which is non-positive for $ t\in[1,2]$ and $ u\in[1,2]$. In fact
$$
\begin{aligned}
 (b+1)\e^{-(b+2)- (t-1)}- (b+u)\e^{-(b+u)(t-1)} &\le  (b+1)\e^{-(b+2)- (1-1)}- (b+u)\e^{-(b+u)(2-1)}\\
 &= (b+1)\e^{-(b+2)}- (b+u)\e^{-(b+u)} \le 0.
\end{aligned}
$$
Hence 
for fixed $u\in[1,2]$, the function in \eqref{thefun} is   non-increasing for $t\in[1,2]$. For  $1 \le x < y \le 2$, we have
\begin{equation}\label{exprel}
 \e^{-(b+u)(x-1)} - \e^{-(b+u)(y-1)} \ge (b+1) \e^{-(b+2)}\left ( \e^{-(x-1)} - \e^{-(y-1)}  \right ). 
\end{equation}
We use this inequality to get  a lower bound for the probability of the   event $\{ Z_2  \in  (x,y), \;  Z_1 \in I_{\eps}(z) \}$. Our strategy is to calculate the probability of  a suitable subset of the latter set. Consider the following event. Suppose that 
\begin{itemize}
\item[a)] $T_{1}< 1$, then
\item[b)] the process spends at $X_{T_{1}}$  an amount of time enclosed   in $(z-1-\eps, z-1+ \eps)$, then
\item[c)] it  jumps to a vertex at level $2$, spends there an amount of time $t$ where  $t+1  \in (x,y)$, and 
\item[d)] it  jumps to level 3 and never returns to  $X_{T_{2}}$.
 \end{itemize}
  In the event just described, levels 1 and 2 are the first two cut levels, and $ \{ Z_2  \in  (x,y) , \;  Z_1 \in I_{\eps}(z) \}$ holds. The probability that  a) holds is exactly $\e^{-b}$. Given $T_{1}=s-1$, the time spent in $X_{T_{1}}$ before the first jump is exponential with parameter $(b+s)$.  Hence b) occurs with probability larger than 
 $$
 \inf_{s \in [1,2]} \Big( \e^{-(b+s)(z - \eps) } - \e^{-(b+s)(z + \eps)} \Big).   
 $$ 
 Given a) and b), the process jumps to level 2 and then to level 3 with probability larger than $ \big( b/(b+2)\big) \big(b/(b+z+\eps))$. The conditional probability, given a) and b),  that the time gap  between these two jumps lies in $(x-1, y-1)$ is larger than
$$ 
  \inf_{u \in  I_{\eps}(z)} \Big( \e^{-(b+u)(x-1) } - \e^{-(b+u)(y-1)} \Big).
$$
At this point, a lower bound for the conditional probability that the process never returns  to $ X_{T_{2}}$ is 
$$
\frac b{b+y}  (1- \alpha_{b}) \ge \frac b{b+2}  (1- \alpha_{b}).
$$
We have
\begin{equation}\label{balletto1}
\begin{aligned}
&\mathbb{P} \Big( Z_2  \in  (x,y), \;  Z_1 \in I_{\eps}(z)  \Big)  \\
&\ge  \e^{-b}\frac {b^{3}}{(b+2)^{2} (b+z+\eps) } \inf_{s \in [1,2]}   \Big( \e^{-(b+s)(z - \eps) } - \e^{-(b+s)(z + \eps)} \Big) \\
&\qquad \qquad \qquad   \inf_{u \in  I_{\eps}(z)} \Big( \e^{-(b+u)(x-1) } - \e^{-(b+u)(y-1)} \Big)(1-\alpha_{b})\\
&\ge (1-\alpha_{b}) \e^{-b}\frac {b^{3} (b+1)}{(b+2)^{2} (b+z+\eps) }   \e^{-(b+2)}\left ( \e^{-(x-1)} - \e^{-(y-1)}  \right ) \inf_{s \in [1,2]} \Big( \e^{-(b+s)(z - \eps) } - \e^{-(b+s)(z + \eps)} \Big),
\end{aligned}
\end{equation}
where in the last inequality we used \eqref{exprel}. Notice that there exists a constant $C_b^{\ssup 4}>0$ such that 
\begin{equation}\label{balletto2}
\inf_{\eps \in (0,1)}   \inf_{z,s \in [1,2]} \frac 1{\eps}\Big( \e^{-(b+s)(z - \eps) } - \e^{-(b+s)(z + \eps)}\Big) \ge C_b^{\ssup 4}.
\end{equation} 
Summarizing, we have 
\begin{equation}\label{balletto3}
\begin{aligned}
\mathbb{P} \Big( Z_2  \in  (x,y), \;  Z_1 \in I_{\eps}(z)  \Big) \ge C_b^{\ssup 5} \left ( \e^{-(x-1)} - \e^{-(y-1)}  \right ) \eps,
\end{aligned}
\end{equation}
where $C_b^{\ssup 5}$ depends only on $b$.

In order to find a lower bound for \eqref{crib} we need to prove that
\begin{equation}\label{balletto4}
\sup_{\eps \in (0,1)}\frac 1\eps\mathbb{P} \Big(   Z_1 \in I_{\eps}(z)  \Big) \le C_b^{\ssup 6},
\end{equation}
for some positive constant $C_b^{\ssup 6}$. To see this, recall the definition of $ B_{j}$ from the proof of Theorem~\ref{codaelle}, and $\zeta$ from \eqref{offwthfc}. The event that level $i$ is not a cut level is subset of $(B_{i}\cap A_{i})^{c}$ (see the proof of Theorem~\ref{codaelle}). Denote by $ m_{i} = h_{1}\big(X_{T_{i}},{\rm fc}(X_{T_{i}})\big)$, which is exponential with mean $1/b$.  Then
$$
\begin{aligned}
\mathbb{P}(Z_{1} \in I_{\eps}(z)) &\le\sum_{i}^{\infty}  \mathbb{P}\Big(m_{i} \in  I_{\eps}(z)\Big) \mathbb{P}\Big(\bigcap_{k=1}^{i-1} (B_{k}\cap A_{k})^{c}\;|\; m_{i} \in I_{\eps}(z)\Big)\\
&\le C \eps \sum_{i}^{\infty} \mathbb{P}\Big(\bigcap_{k=1}^{i-1} (B_{k}\cap A_{k})^{c}\;|\; m_{i} \in I_{\eps}(z)\Big),
\end{aligned}
$$
where the constant $C$ is independent of $\eps$ and $z$. It remains to prove that the sum in the right-hand side is bounded by a constant independent of $\eps$. Notice that, for $ i > \zeta$, $ A_{i - \zeta}$ and  $ B_{i - \zeta}$ are independent of $m_{i}$. Moreover the events
$$ 
A_{i - \zeta}\cap B_{i - \zeta}, \; A_{i - 2\zeta}\cap B_{i - 2\zeta}, \; A_{i - 3\zeta}\cap B_{i -3\zeta},\ldots 
$$
are independent by the proof of Proposition~\ref{aind}. Hence
$$
\begin{aligned}
\mathbb{P}(Z_{1} \in I_{\eps}(z)) &\le C \eps \sum_{i}^{\infty} \mathbb{P}\Big(\bigcap_{k=1}^{[(i-1)/\zeta]} (B_{i-k\zeta}\cap A_{i - k\zeta})^{c}\;|\; m_{i} \in I_{\eps}(z)\Big)\\
&=  C \eps \sum_{i}^{\infty} \mathbb{P}\Big(\bigcap_{k=1}^{[(i-1)/\zeta]} (B_{i-k\zeta}\cap A_{i - k\zeta})^{c}\Big) \qquad \mbox{(by independence)}\\
&= C \eps \sum_{i}^{\infty} \mathbb{P}\Big((B_{i-k\zeta}\cap A_{i - k\zeta})^{c}\Big)^{[(i-1)/\zeta]}<\infty.
\end{aligned}
$$
Combining \eqref{balletto1},  \eqref{balletto2} and \eqref{balletto4}, we get
\begin{equation}\label{ekk}
\begin{aligned}
\mathbb{P} \Big( Z_2 \; \in \; (x,y) \; \big{|} \;  Z_1 =z   \Big) &= \lim_{\eps \downarrow 0}\frac 1{\mathbb{P} \Big(   Z_1 \in I_{\eps}(z)  \Big)}\mathbb{P} \Big( Z_2 \; \in \; (x,y) \; , \;  Z_1 \in I_{\eps}(z)  \Big)\\
 &\ge \lambda  \frac{ \left ( \e^{-(x-1)} - \e^{-(y-1)}  \right )}{\left( 1 - \e^{-1} \right)},
\end{aligned}
\end{equation}  
\noindent for some $\lambda >0$.
. A finite measure defined on  field $\mathcal{A}$ can be extended uniquely to the sigma-field  generated by $\mathcal{A}$, and this extension coincides with the outer measure. We apply this result to prove that \eqref{ekk}  holds for any Borel set $C \subset [1,2]$, using the fact that it holds  in the field of finite unions of intervals.
For any interval $E$, the right-hand side of \eqref{ekk} can be written in an integral form as
 $$\lambda  \int_{E}  \frac{ e^{-x+1} }{\left( 1 - \e^{-1} \right)} \d x.$$
Fix a Borel set $C\subset[1,2] $ and  $\eps >0$ choose a countable collection of disjoint intervals $E_{i}\subset[1,2]$, $ i \ge 1$, with $ C \subset \bigcup_{i=1}^{\infty} E_{i}$, such that 
$$
\begin{aligned}
 \mathbb{P} (Z_2 \in  C\; \big{|} \;  Z_1 =z  ) &\ge  \sum_{i=1}^{\infty} \mathbb{P} (Z_2 \in  E_{i} \; \big{|} \;  Z_1 =z  ) - \eps \\
 &\ge \lambda \sum_{i=1}^{\infty}  \int_{E_{i}}  \frac{ e^{-x+1} }{\left( 1 - \e^{-1} \right)} \d x - \eps\\
 &\ge \lambda \int_{C}  e^{-x+1}/\left( 1 - \e^{-1} \right) \d x - \eps.
\end{aligned}
$$
The first inequality is true because of the extension theorem, and the fact that the right-hand side is a lower bound for the outer measure, for a suitable choice of the $E_{i}$s.
 The inequality \eqref{homzeta}, with $\phi(C) = \int_{C}  e^{-x+1}/\left( 1 - \e^{-1} \right) \d x$,   follows by
sending $\eps$ to $0$.
\qed           \end{proofsect}
The proof of the following Proposition can be found in \cite{AN1978}.
\begin{proposition}
There exists a constant $\varrho \in (0,1)$ and a sequence  of random times 
$\{N_k, \, k \ge 0\}$, with $N_{0}=0$,   such that 
\begin{itemize}
\item the sequence  $\{ Z_{N_k},\, k \ge 1\} $ consists of   independent and identically distributed random variables with distribution $\phi(\cdot)$
\item $N_i- N_{i-1}$,  $i \ge 1$, are i.i.d.  with a geometric distribution($\rho$), i.e.
$$ \mathbb{P}(N_{2} - N_{1} = j) = (1 - \varrho)^{j-1} \varrho, \qquad \mbox{ with $j \ge 1$}. $$
\end{itemize}
\end{proposition}
\begin{Lemma} $ \sup_{i \in \N} \mathbb{E} [ (\tau_{N_{i+1}} - \tau_{N_{i}})^{2} ] < \infty$.
\end{Lemma}
\begin{proofsect}{Proof} 
It is enough to prove  $\mathbb{E} [ (\tau_{N_2} - \tau_{N_{1}})^{2} ] < \infty $.
By virtue of Jensen's inequality, we have that  
\begin{equation}
\begin{aligned}
\mathbb{E}\Big[ (\tau_{k}-\tau_{m})^{11/5}\Big] &= \mathbb{E}\Big[ (\sum_{j=1}^{k-m} \tau_{m+j}-\tau_{m+j-1})^{11/5}\Big]\\
&\le (k-m)^{11/5 } \mathbb{E}[(\tau_{2} -\tau_{1})^{11/5}].
\end{aligned}
\end{equation}
 Using Holder with $p = 11/10$, we have
\begin{equation*}
\begin{aligned}
& \mathbb{E} [ (\tau_{N_2} - \tau_{N_{1}})^{2} ]  = \sum_{k=2}^{\infty} \sum_{m=1}^{k-1} \mathbb{E} \left [ (\tau_k - \tau_{m})^{2} \1_{\{ N_1= m,\, N_{2} = k  \}} \right] \\
&\le  \sum_{k=2}^{\infty} \sum_{m=1}^{k-1}   \mathbb{E} \Big [  \bigl(\tau_{k} - \tau_{m} \bigr )^{11/5} \Big]^{ 10/11} \mathbb{P} ( N_1 = m,\, N_{2} - N_{1}= k - m )^{1/11 } \qquad  \\
&= \sum_{k=2}^{\infty} \sum_{m=1}^{k-1}   \mathbb{E} \Big [  \bigl(\tau_{k} - \tau_{m} \bigr )^{11/5} \Big]^{ 10/11}  \mathbb{P} (N_1 = m)^{1/11 }  \mathbb{P} (N_{2} - N_{1}= k - m )^{1/11 }  \\
&\le \sum_{k=2}^{\infty}  \sum_{m=1}^{k-1} (k-m)^{3} \mathbb{E}[(\tau_{2} -\tau_{1})^{11/5}]^{ 10/11} \varrho^{2/11} \big(1-\varrho\big) ^{(k-2)/11 }\\
&\le \varrho^{2/11} \mathbb{E}[(\tau_{2} -\tau_{1})^{11/5}]^{ 10/11}  \sum_{k=2}^{\infty}  k^{4} \big(1-\varrho\big) ^{(k-2)/11 } 
< \infty, 
\end{aligned}
\end{equation*}
 where we used the fact that $ 0<  \varrho < 1$. 
\qed           \end{proofsect}
With a similar proof we get the following result.
\begin{Lemma} $ \sup_{i \in \N} \mathbb{E} \big[\big( l_{N_{i+1}} - l_{N_{i}} \big)^{2}\big] < \infty.$
\end{Lemma}
\begin{definition}
 A process $\{Y_k,\, k \ge 1\}$, is said to be {\bfseries one-dependent} if $Y_{i+2}$ is independent of $\{ Y_j, \mbox{ with } 1 \le  j \le i\}$. 
\end{definition}
\begin{Lemma}
Let $\Upsilon_{i} := \big(\tau_{N_{i+1}} - \tau_{N_{i}},\, l_{N_{i+1}}- l_{N_{i}}\big) $, for $i \ge 1$. The process  $ \mathbf{\Upsilon}:= \big\{ \Upsilon_{i},\, i \ge 1 \big\} $ is one-dependent. Moreover $ \Upsilon_{i}$, $i \ge 1$, are    identically distributed.  
\end{Lemma}
\begin{proofsect}{Proof} 
Given $Z_{N_{i-1}}$, $\Upsilon_{i} $ is  independent of $\{ \Upsilon_{j} , \, j \le i-2\}$.  Thus,  it is sufficient to prove that 
 $\Upsilon_{i} $ is independent of $Z_{N_{i-1}}$. To see this,  it is enough 
to realize that given $Z_{N_{i}}$, $\Upsilon_{i} $ is independent of $Z_{N_{i-1}}$, and
 combine this with the fact that $Z_{N_{i}}$ and $Z_{N_{i-1}}$ are independent.   The variables $ Z_{N_i}$ are i.i.d., hence $\{\Upsilon_{i}, \, i \ge 2\}$, are identically distributed. \hspace{\fill} 
\qed           \end{proofsect}

  The Strong Law of Large Numbers holds for one-dependent  sequences of identically distributed variables bounded in $\Lcal^1$. To see this, just consider separately the sequence of random variables with even and odd indices and apply the usual Strong Law of Large Numbers to each of them. 

\noindent  Hence, for some constants  $0 < C_b^{\ssup 7}, C_b^{\ssup 8} < \infty$, we have
\begin{equation}\label{stronglaws}
 \lim_{i \rightarrow \infty} \frac{\tau_{N_i}}{i} \rightarrow C_b^{\ssup 7},\qquad \mbox{and} \qquad  \lim_{i \rightarrow \infty} \frac{l_{N_i}}{i} \rightarrow C_b^{\ssup 8}, \qquad a.s.. 
 \end{equation}
\begin{proofsect}{Proof of Theorem 1}
 If $\tau_{N_i} \le t < \tau_{N_{i+1}}$, then by the definition of cut level, we have
\[ l_{N_i} \le |X_t| < l_{N_{i+1}}. \]
Hence
\[  \frac{l_{N_i}}{\tau_{N_{i+1}}} \le \frac{|X_t|}{t} < \frac{l_{N_{i+1}}}{\tau_{N_{i}}}.\]
\noindent Let 
\begin{equation}\label{constK}
 K_{b}^{\ssup 1} = \frac{\mathbb{E}[l_{N_{2}} - l_{N_{1}}]}{\mathbb{E}[\tau_{N_{2}} - \tau_{N_{1}}]},
\end{equation}
  which are the constants in \eqref{stronglaws}. Then 
\[ \limsup_{t \rightarrow \infty} \frac{ |X_t|}{t} \le \lim_{i \rightarrow \infty} \frac{l_{N_{i+1}}}{ \tau_{N_i}}  = \lim_{i \rightarrow \infty} \frac{l_{N_{i+1}}}{ i+1} \frac{i}{\tau_{N_i}} = K^{\ssup 1}_{b}, \; \mbox{a.s..}  \]
\noindent Similarly, we can prove that
\[ \liminf_{t \rightarrow \infty} \frac{ |X_t|}{t} \ge K^{\ssup 1}_{b}, \;  \mbox{a.s..} \]  
Now we turn to the proof of the central limit theorem. 
First we prove that there exists a constant $C\ge 0$ such that
\begin{equation}\label{cltforl}
 \frac{l_{N_{m}} -  K_{b}^{\ssup 1}\tau_{N_{m}}}{\sqrt m} \Longrightarrow  \mbox{Normal$(0,C)$},
\end{equation}
where Normal$(0,0)$ stands for  the  Dirac mass at 0.
To prove \eqref{cltforl} we use a theorem  from  \cite{HR1948}. The reader can find the  statement of this theorem in the Appendix, Theorem \ref{Hoef-Rob}, (see also \cite{Ser1968}). In order to apply this result we first need to prove that the quantity
\begin{equation}\label{econv}
  \frac 1m \mathbb{E}\Big[ \big( l_{N_{m}} - K^{\ssup 1}_{b}\tau_{N_{m}} \big)^{2}\Big] = \mathbb{E}\Big[ \Big( \frac{l_{N_{m}} - K^{\ssup 1}_{b}\tau_{N_{m}}}{\sqrt{m}}\Big)^{2} \Big]
\end{equation}
converges.  Call $Y_{1} =l_{N_{1}}-K^{\ssup 1}_{b}(\tau_{N_{1}}$ and let  $Y_{i} = l_{N_{i}}- l_{N_{i-1}}-K^{\ssup 1}_{b}(\tau_{N_{i}}-\tau_{N_{i-1}})$, with $ i \ge 2$. The quantity in \eqref{econv}
can be written as
$$ \frac 1m \mathbb{E}\Big[ \big( \sum_{i=1}^{m}Y_{i} \big)^{2}\Big].$$ 
The random variables $Y_{i}$ are identically distributed with the exception of $Y_{1}$. From the definition of $ K^{\ssup 1}_{b}$ given in \eqref{constK}, we have
$$\mathbb{E}[Y_{i}] = \mathbb{E}[l_{N_{2}} - l_{N_{1}}] - \mathbb{E}[l_{N_{2}} - l_{N_{1}}]=0.$$

Hence $Y_{i}$, with $i\ge1$,  is a zero-mean  one-dependent process, and we get
\begin{equation}
\begin{aligned}
\mathbb{E}\Big[ &\big( l_{N_{m}} - K^{\ssup 1}_{b}\tau_{N_{m}} \big)^{2} \Big]= \mathbb{E}\Big[ \Big(\sum_{i=1}^{m} Y_{i} \Big)^{2} \Big] \\
&= (m-1) \mathbb{E}[Y_{2}^{2}] 
+ 2(m-2)\mathbb{E}[Y_{3}Y_{2 }] + \mathbb{E}[Y_{1}^{2}] + 2 \mathbb{E}[Y_{1}Y_{2}].
\end{aligned}
\end{equation}
This proves that the limit in \eqref{econv} exists and is equal to $\mathbb{E}[ Y_{2}^{2}] + 2\mathbb{E}[Y_{3}Y_{2 }]$. 
 Now we face two options. If the limit is equal to zero, then using Chebishev we get that
$$
\lim_{m \ti}\mathbb{P}\Big( \Big| \frac{l_{N_{m}} - C\tau_{N_{m}}}{\sqrt{m}} \Big|> \eps\Big)=\lim_{m \ti}\mathbb{P}\Big( \Big| \frac1{\sqrt{m}} \sum_{i=1}^{m} Y_{i} \Big|> \eps\Big)  \le  \lim_{m \ti} \frac 1{\eps}\mathbb{E}\Big[ \Big( \frac{\sum_{i=1}^{m}Y_{i}}{\sqrt{m}}\Big)^{2} \Big]=0.$$
If the limit of the quantity in \eqref{econv} is positive, then we can apply Theorem~~\ref{Hoef-Rob} and deduce central limit theorem for $Y_{i}, i \ge 1$, yielding \eqref{cltforl}. 

Now we use \eqref{cltforl} to prove the central limit theorem for $|X_{t}|$.
 If $\tau_{N_m} \le t < \tau_{N_{m+1}}$, then
\begin{equation}\label{finalclt}
\begin{aligned}
\frac {|X_{t}| - K_{b}^{\ssup 1}t}{K^{\ssup 2}_{b}\sqrt{t}} &\ge \frac{l_{N_{m}} - K_{b}^{\ssup 1}\tau_{N_{m+1}}}{K^{\ssup 2}_{b}\sqrt{\tau_{N_{m+1}}}}= \sqrt{\frac m{\tau_{N_{m+1}}}}\Big(\frac{l_{N_{m}}- K_{b}^{\ssup 1}\tau_{N_{m}}}{\sqrt{m}} + \frac{K^{\ssup 1}_{b}}{\sqrt{m}}(\tau_{N_{m}}- \tau_{N_{m+1}})\Big)\\
&= \sqrt{\frac m{\tau_{N_{m+1}}}}\Big(\frac{\sum_{i=1}^{m}Y_{i}}{\sqrt{m}} - \frac{Y_{m}K^{b}_{1}}{\sqrt{m}}\Big).
\end{aligned}
\end{equation}
The last expression converges, by virtue of the Slutzky's lemma, either to a Normal distribution or to a Dirac mass at 0, depending on whether  the limit in \eqref{econv} is positive or is zero. To see this, notice that
\begin{eqnarray*}
\lim_{m \ti}\sqrt{\frac m{\tau_{N_{m+1}}}}&=&\sqrt{\frac{1}{\mathbb{E}[\tau_{N_{2}}- \tau_{N_{1}}]}}, \qquad \mbox{a.s.}\\ 
\frac{\sum_{i=1}^{m}Y_{i}}{\sqrt{m}} &\Longrightarrow& \mbox{Normal$(0,C)$}\\
\lim_{m \ti} \frac{Y_{m}K^{b}_{1}}{\sqrt{m}}&=&0, \qquad \mbox{a.s.}. 
\end{eqnarray*}
Similarly
$$ \frac {|X_{t}| - K_{b}^{\ssup 1}t}{K^{\ssup 2}_{b}\sqrt{t}} \le \sqrt{\frac {m+1}{\tau_{N_{m}}}}\Big(\frac{\sum_{i=1}^{m+1}Y_{i}}{\sqrt{m+1}} + \frac{Y_{m+1}K^{b}_{1}}{\sqrt{m}}\Big),$$
and the right-hand side converges to the same limit of the right-hand side of \eqref{finalclt}.\qed
\end{proofsect}
 \section{Appendix} We include a corollary to a result of Hoeffding and Robbins (see \cite{HR1948} or \cite{Ser1968}).
\begin{theorem}[Hoeffding-Robbins]\label{Hoef-Rob} Suppose $\mathbf{Y}:= \{Y_{i},\, i \ge 1\}$ is a one-dependent process whose components are identically distributed with mean 0. If 
\begin{itemize}
\item $\mathbb{E}[ Y_{i}^{2+\delta}] < \infty,$     for some  $\delta >0,$
\item $\lim_{n \ti }\frac1n {\rm Var}(\sum_{i=1}^{n} Y_{i}) $ converges to a positive finite constant $K$, then
\end{itemize} 
$$ \frac{\sum_{i=1}^{n} Y_{i} - n\mathbb{E}[Y_{1}]}{K \sqrt{n}} \Longrightarrow {\rm Normal}(0,1).$$
\end{theorem}

\noindent{\bf Aknowledgement.} The author was supported by the DFG-Forschergruppe 718  ''Analysis and stochastics in complex physical systems'', and by the Italian PRIN 2007 grant 2007TKLTSR "Computational
 markets design and agent-based models of trading behavior".  The author would like to thank Burgess Davis and an anonymous referee for helpful suggestions.


\begin{thebibliography}{}
%
%
\bibitem{A2008}
E. Aid\'ekon. Transient random walks in random environment on a Galton-Watson tree. \textit{Prob. Th. Rel. Fields.} \textbf{142,} (2008) 525-559.
\bibitem{AN1978}
K. Athreya, P. Ney. A new approach to the limit theory of recurrent markov chains.\textit{Trans. American Math. Society}   \textbf{245,} (1978) 493-501.
\bibitem{C2006}
A. Collevecchio. Limit theorems for reinforced random walks on   certain trees. \textit{Prob. Th. Rel. Fields.} \textbf{136} (2006) 81-101.
\bibitem{C2006b}
A. Collevecchio. On the transience of processes defined on Galton-Watson trees. \textit{ Ann. Probab.} \textbf{34} (2006) 870-878.
\bibitem{C1987}
D. Coppersmith, P. Diaconis. Random walks with reinforcement. \textit{Unpublished manuscript} (1987).
\bibitem{D1990}
B. Davis. Reinforced random walk.\textit{ Prob. Th. Rel. Fields}   \textbf{84,} (1990) 203-229.
\bibitem{D1999}
B. Davis. Reinforced and perturbed random walks. In P. R\'ev\'esz and B. T\'oth (Eds.) \textit{Random walks}, Volume   \textbf{9,}  113-126 (Bolyai Soc. Math. Studies).
\bibitem{DV2002}
B. Davis, S. Volkov. Continuous time vertex-reinforced jump processes.\textit{ Prob Th. Rel. Fields},   \textbf{84,} (2002) 281-300.
\bibitem{DV2004}
B. Davis, S. Volkov. Vertex-reinforced jump process on trees and finite graphs.  \textit{Prob. Th. Rel. Fields}   \textbf{128,} (2004) 42-62.
\bibitem{DZ}
A.Dembo, O. Zeitoni. \textit{Large deviations techniques and applications} (Springer, 1998). 
\bibitem{DKL2002}
R. Durrett , H. Kesten, V. Limic. Once reinforced random walk. \textit{ Prob. Th. Rel. Fields}   \textbf{122,} (2002) 567-592.
\bibitem{LT2007}
V. Limic, P. Tarr\`es. Attracting edge and strongly edge reinforced walks. \textit{Ann. Probab.} \textbf{35} (2007) 1783-1806.
\bibitem{MR2008}
F. Merkl, S.W.W. Rolles. Recurrence of edge-reinforced random walk on a two-dimensional graph. \textit{Preprint} (2008).
\bibitem{MPU2006}
F. Merlev\`ede, M. Peligrad, S. Utev. Recent advances in invariance principles for stationary sequences. \textit{Probability Surveys} \textbf{3} (2006) 1-36.
\bibitem{HR1948}
W. Hoeffding, H. Robbins. The central limit theorem for dependent random variables. \textit{Duke. Math. J.}  \textbf{15} (1948)
773-780.
\bibitem{P1988}
R. Pemantle. Phase transition in reinforced random walks and rwre on trees. \textit{Ann. Probab.}   \textbf{16,} (1988) 1229-1241.
\bibitem{P1992}
R. Pemantle. Vertex-reinforced random walk. \textit{Prob. Th. Rel. Fields} \textbf{92,} (1992) 117-136.
\bibitem{P2007}
R. Pemantle. A survey of random processes with reinforcement. \textit{Probability Surveys} \textbf{4} (2007) 1-79.
\bibitem{PV99}
R. Pemantle, S. Volkov. Vertex-reinforced random walk on Z 
has finite range.  \textit{Ann. Probab.} \textbf{27,} (1999) 1368Ð1388.
\bibitem{S1994}
T. Sellke. Reinforced random walk on the d-dimensional integer lattice. \textit{Markov Process. Relat. Fields}   \textbf{14,} (2008) 291-308.
\bibitem{S2006}
T. Sellke. Recurrence of reinforced random walk on a ladder. \textit{Electr. Journal of Probab.} \textbf{11,} (2006) 301-310.
\bibitem{Ser1968}
R. J. Serfling. Contributions to the central limit theory for dependent variables. \textit{Ann. of Math. Statistics} \textbf{39,}(1968) 1158-1175.
\bibitem{T2004}
P. Tarr\`es. Vertex-reinforced random walk on Z eventually gets stuck on five points.   \textit{Ann. Probab.} \textbf{32,}  (2004) 2650-2701.
\bibitem{V2001}
S. Volkov. Vertex-reinforced random walk on arbitrary graphs. \textit{Ann. Probab.} \textbf{29,}, (2001) 66-91.
\end{thebibliography}
\end{document}